\pdfoutput=1
\documentclass[12pt]{amsart}

\usepackage{amssymb,mathtools}
\usepackage{booktabs,array}
\usepackage{enumitem}
\usepackage[hidelinks]{hyperref}
\hypersetup{pdfauthor={Hyunjun Yi},pdftitle={Isosceles sets and two-distance sets in Euclidean space}}

\makeatletter
\@namedef{subjclassname@2020}{%
  \textup{2020} Mathematics Subject Classification}
\makeatother

\newtheorem{theorem}{Theorem}[section]
\newtheorem{lemma}[theorem]{Lemma}
\newtheorem{proposition}[theorem]{Proposition}
\newtheorem{corollary}[theorem]{Corollary}
\newtheorem{problem}[theorem]{Problem}
\theoremstyle{definition}

\newtheorem{remark}[theorem]{Remark}
\newtheorem{example}[theorem]{Example}

\numberwithin{equation}{section}

\usepackage[letterpaper,textwidth=13.5cm,textheight=23cm,centering]{geometry}
\newcommand{\R}{\mathbb{R}}

\newcommand{\ip}[2]{\langle #1,#2\rangle}
\newcommand{\nr}[1]{\lVert #1\rVert^2}
\DeclareMathOperator{\aff}{aff}
\DeclareMathOperator{\iso}{iso}
\DeclareMathOperator{\pyr}{pyr}
\DeclareMathOperator{\tr}{tr}
\newcommand{\J}{\mathcal{J}}
\newcommand{\LL}{\mathcal{L}}
\newcommand{\B}{\mathcal{B}}

\begin{document}

\title[Isosceles sets and two-distance sets]{Isosceles sets and two-distance sets in Euclidean space}

\author[H. Yi]{Hyunjun Yi}
\address{Seoul International School, 15 Seongnam-daero 1518beon-gil, Sujeong-gu, Seongnam-si, Gyeonggi-do 13113, Republic of Korea}
\email{hjyi2027@gmail.com}

\date{}

\begin{abstract}
A finite set of points in $\R^d$ is isosceles if every three of its points form an isosceles triangle. Erd\H{o}s asked in 1946 for the largest size $\iso(d)$ of an isosceles set in $\R^d$; the answer was known for $d\le8$. Let $g(d)$ and $M(d)$ be the largest sizes of a two-distance set and of a spherical two-distance set in $\R^d$. We prove that $\iso(d)=\max\{g(d),M(d)+1,M(d-1)+3\}$ for every $d\ge2$, which reduces the problem to these two classical ones. With the recently determined values of $g(d)$ for $9\le d\le22$ we obtain $\iso(d)=\binom{d+1}{2}+1$ for $9\le d\le21$ and $\iso(22)=276$, and we show that the extremal set is unique: for $d\le21$ it is the set of midpoints of the edges of a regular simplex together with its centroid, and for $d=22$ it is the McLaughlin configuration together with its centre. For $10\le d\le22$ every isosceles set with more than $\binom d2+3$ points is a two-distance set or a spherical two-distance set together with its centre. The proofs start from Blokhuis' decomposition of isosceles sets into mutually orthogonal two-distance sets. Counting leaves one configuration open: a two-distance set with $\binom{d+1}{2}$ points in a hyperplane, the largest number allowed by Blokhuis' bound, together with two points that are mirror images of each other. We exclude it by showing that a two-distance set $T$ attaining Blokhuis' bound in $\R^n$ is unisolvent for quadratic polynomials and that, if $T$ lies on two concentric spheres, then apart from two exceptional cases its two layers are spherical $2$-designs whose sizes are determined by the radii. We also show that the squared distances of such a set are in ratio $(k-1):k$ with $(2k-1)^2\le n+1$ if $n\ge9$. Hence there is none for $9\le n\le23$, so that $g(23)\le299$ and $\iso(23)=\max\{g(23),278\}\le299$.
\end{abstract}



\maketitle

\section{Introduction}\label{sec:intro}

A finite set $S$ in a metric space is \emph{isosceles} if among any three points of $S$ at least two of the three pairwise distances coincide. Degenerate triangles, such as three equally spaced points on a line, are allowed. In 1946 Erd\H{o}s \cite{Erd46} asked for the largest cardinality of an isosceles set in the plane and in space, noting that the plane admits six points. Kelly \cite{EK47} proved that the six-point planar set is unique (the vertices and the centre of a regular pentagon) and gave an eight-point isosceles set in $\R^3$, which Croft \cite{Cro62} showed to be maximal and Kido \cite{Kid06} showed to be unique. Write $\iso(d)$ for the largest size of an isosceles set in $\R^d$. Ionin \cite{Ion09} determined $\iso(d)$ for $d\le7$, described all maximal isosceles sets in these dimensions, and observed that $\iso(8)=45$ (Kido \cite{Kid10} also determined $\iso(4)=11$ and showed that there are exactly two extremal sets), so that
\[
\iso(1),\dots,\iso(8)=3,\ 6,\ 8,\ 11,\ 17,\ 28,\ 30,\ 45
\]
(sequence A175769 of \cite{OEIS}; the problem is Erd\H{o}s Problem \#503 in \cite{EP503}). No value $\iso(d)$ with $d\ge9$ was known.

Isosceles sets are closely related to \emph{two-distance sets}, sets whose pairwise distances take only two values: every two-distance set is isosceles, and Blokhuis \cite[Chapter~7]{Blo83} showed that every isosceles set is built from two-distance sets lying in mutually orthogonal subspaces (Section~\ref{sec:decomp}). Let $g(d)$ be the largest size of a two-distance set in $\R^d$. The values $g(d)$ were determined for $d\le8$ by Einhorn and Schoenberg \cite{ES66} and Lison\v{e}k \cite{Lis97}, and for $9\le d\le22$ in the recent paper \cite{BLSTYY}:
\begin{equation}\label{eq:g}
\begin{gathered}
g(1),\dots,g(9)=3,5,6,10,16,27,29,45,45,\\
g(d)=\binom{d+1}{2}\quad(10\le d\le21),\qquad g(22)=275.
\end{gathered}
\end{equation}
For $10\le d\le21$ the value $\binom{d+1}{2}$ is attained, up to similarity, only by the set
\[
\J(d+1,2)=\{e_i+e_j:1\le i<j\le d+1\}\subset\R^{d+1},
\]
the midpoints of the edges of a regular simplex. In $\R^9$ the set $\J(10,2)$ is the only $45$-point two-distance set that spans $\R^9$, while Lison\v{e}k's $45$-point set $\LL$ spans an $8$-dimensional subspace. Since $\J(d+1,2)$ lies on a sphere, adding its centroid gives an isosceles set with $\binom{d+1}{2}+1$ points in $\R^d$. In $\R^{22}$ the McLaughlin configuration is a spherical two-distance set with $275$ points \cite{DGS77}, and adding its centre gives an isosceles set with $276$ points. Our main result is that nothing better is possible.

\begin{theorem}\label{thm:main}
\begin{enumerate}[label=\textup{(\alph*)}]
\item For $9\le d\le 21$, $\iso(d)=\binom{d+1}{2}+1$. The unique isosceles set of this size in $\R^d$, up to similarity, is $\J(d+1,2)$ together with its centroid.
\item $\iso(22)=276$. The unique isosceles set of this size in $\R^{22}$, up to similarity, is the McLaughlin configuration together with the centre of its circumsphere.
\end{enumerate}
\end{theorem}

Theorem~\ref{thm:main} extends the sequence A175769 by fourteen terms:
\begin{gather*}
\iso(9),\dots,\iso(15)=46,\,56,\,67,\,79,\,92,\,106,\,121,\\
\iso(16),\dots,\iso(22)=137,\,154,\,172,\,191,\,211,\,232,\,276 .
\end{gather*}
The upper bound is a consequence of Blokhuis' decomposition, of \eqref{eq:g}, and of the known maximum sizes of spherical two-distance sets: $d(d+1)/2$ for $7\le d\le21$ (Musin \cite{Mus09}) and $275$ for $d=22$ \cite{DGS77}. For $10\le d\le22$ these ingredients give the exact value directly. In dimension $9$ they give only $\iso(9)\le47$, the value $47$ coming from a hypothetical configuration consisting of two points together with a $45$-point two-distance set, such as Lison\v{e}k's set, in an orthogonal $8$-space. We exclude it in Section~\ref{sec:nine} by showing that no isosceles set has a decomposition whose last part is an $8$-dimensional two-distance set with at least $44$ points. This argument uses only general bounds for two-distance sets of dimension at most $8$ together with the results of \cite{BLSTYY} in dimension $9$, and it does not require a classification of the large two-distance sets in $\R^8$.

A note posted on the discussion forum of Erd\H{o}s Problem \#503 \cite{EP503} states the identity
\begin{equation}\label{eq:forum}
\iso(d)=\max\{g(d),\ M(d)+1,\ M(d-1)+3\}\qquad(d\ge2),
\end{equation}
where $M(d)$ is the largest size of a spherical two-distance set in $\R^d$. We have not relied on that note. Combined with \eqref{eq:g} it would give the values in Theorem~\ref{thm:main}. The present paper gives a complete proof of these values from \eqref{eq:g}, including the uniqueness statements and the structure theorem below. The inputs \eqref{eq:g} for $d\ge9$, on which any such evaluation depends, are from \cite{BLSTYY}. In the second part of the paper we give a complete proof of \eqref{eq:forum} in every dimension (Theorem~\ref{thm:all} below).

The same method yields a description of all large isosceles sets.

\begin{theorem}\label{thm:structure}
Let $10\le d\le 22$ and let $S\subset\R^d$ be an isosceles set with more than $\binom d2+3$ points. Then $S$ is a two-distance set, or $S$ is a spherical two-distance set together with the centre of its sphere. The bound is sharp: $\J(d,2)$ together with its centre and the two points on the axis of its sphere at distance equal to its radius from the centre is an isosceles set with $\binom d2+3$ points in $\R^d$ which is of neither type. In $\R^9$ the same conclusion holds for isosceles sets with more than $45$ points.
\end{theorem}

The method also determines $\iso(23)$ in terms of $g(23)$, the first unknown value of $g$. In $\R^{23}$ Ge, Koolen and Munemasa \cite{GKM26} constructed a two-distance set with $277$ points, and Blokhuis' bound \cite{Blo84} gives $g(23)\le\binom{25}{2}=300$; we improve this to $g(23)\le299$ (Corollary~\ref{cor:g23}).

\begin{theorem}\label{thm:23}
$\iso(23)=\max\{g(23),278\}$. In particular $278\le\iso(23)\le299$.
\end{theorem}

The values in Theorems~\ref{thm:main} and~\ref{thm:23} agree with \eqref{eq:forum}. Our second main result is that \eqref{eq:forum} holds in every dimension.

\begin{theorem}\label{thm:all}
For every $d\ge2$,
\[
\iso(d)=\max\{g(d),\ M(d)+1,\ M(d-1)+3\}.
\]
\end{theorem}

Thus the problem of Erd\H{o}s on isosceles sets is reduced, in every dimension, to the two classical problems of determining the largest two-distance sets and the largest spherical two-distance sets. More precisely, an isosceles set in $\R^d$ that is not a two-distance set has at most $\max\{M(d)+1,M(d-1)+3\}$ points (Corollary~\ref{cor:three}). The three terms in Theorem~\ref{thm:all} are attained by a largest two-distance set, by a largest spherical two-distance set together with its centre, and by a largest spherical two-distance set in a hyperplane together with its centre and two further points (Remark~\ref{rem:lower}). By the theorem of Glazyrin and Yu \cite{GY18}, $M(a)=\binom{a+1}{2}$ for all $a\ge7$ that are not of the form $(2k+1)^2-3$, so that Theorem~\ref{thm:all} takes the following form in most dimensions.

\begin{corollary}\label{cor:GY}
Let $d\ge8$, and suppose that neither $d$ nor $d-1$ is of the form $(2k+1)^2-3$ with an integer $k$. Then
\[
\iso(d)=\max\Bigl\{g(d),\ \binom{d+1}{2}+1\Bigr\}.
\]
\end{corollary}

The proof of Theorem~\ref{thm:all} has two parts. Call a two-distance set $X\subset\R^d$ a \emph{two-distance pyramid} with \emph{apex} $p\in X$ if its \emph{base} $X\setminus\{p\}$ is nonempty and lies in a hyperplane $H$ with $p\notin H$, and let $\pyr(d)$ be the largest size of a two-distance pyramid in $\R^d$. Adding to a two-distance pyramid the mirror image of its apex in $H$ gives an isosceles set (Lemma~\ref{lem:double}), so that $\iso(d)\ge\pyr(d)+1$. A two-distance pyramid is a two-distance set, and its base is a two-distance set in a hyperplane, so Blokhuis' bound $g(a)\le\binom{a+2}{2}$ \cite{Blo84} gives
\begin{equation}\label{eq:pyrbounds}
\pyr(d)\le g(d)\qquad\text{and}\qquad\pyr(d)\le g(d-1)+1\le\binom{d+1}{2}+1 .
\end{equation}
The first part of the proof refines the decomposition argument and shows that pyramids are the only obstacle to \eqref{eq:forum}.

\begin{theorem}\label{thm:formula}
For every $d\ge2$,
\[
\iso(d)=\max\{g(d),\ M(d)+1,\ M(d-1)+3,\ \pyr(d)+1\}.
\]
\end{theorem}

Since $M(d)\ge\binom{d+1}{2}$, the last term matters only if $\pyr(d)=\binom{d+1}{2}+1$, that is, only if some two-distance set with $\binom{d+1}{2}$ points in a hyperplane, a set attaining Blokhuis' bound in dimension $d-1$, is the base of a two-distance pyramid. The second part of the proof shows that this is impossible.

\begin{theorem}\label{thm:nopyr}
Let $n\ge1$ and let $T\subset\R^n$ be a two-distance set with $\binom{n+2}{2}$ points and squared distances $a<b$. Then there are no point $m\in\R^n$ and no real number $t$ such that $\nr{x-m}+t\in\{a,b\}$ for all $x\in T$; that is, $T$ does not lie on two concentric spheres whose squared radii differ by $b-a$. Consequently $T$ is not the base of a two-distance pyramid, and $\pyr(d)\le\binom{d+1}{2}$ for every $d\ge2$.
\end{theorem}

Theorem~\ref{thm:all} follows from Theorems~\ref{thm:formula} and~\ref{thm:nopyr}. Apart from three equally spaced points on a line, Lison\v{e}k's set $\LL$ in $\R^8$ is, to our knowledge, the only known two-distance set that attains Blokhuis' bound, and it is not known whether there are others; Theorem~\ref{thm:nopyr} does not depend on the answer. Its proof rests on two properties of a two-distance set $T\subset\R^n$ with $\binom{n+2}{2}$ points, which may be of independent interest. First, every real function on $T$ is the restriction of a unique polynomial of degree at most $2$ (Lemma~\ref{lem:interp}; the argument is also a short proof of Blokhuis' bound). Second, if $T$ lies on two concentric spheres, then, apart from two exceptional cases, the points of $T$ on each of the spheres form a spherical $2$-design, their number is determined by the radius, and the two radii satisfy an explicit equation (Theorem~\ref{thm:twosphere}); the set $\LL$, with $36$ and $9$ points on two spheres, is an example. For the base of a two-distance pyramid these conditions contradict each other.

We also restrict the two-distance sets that attain Blokhuis' bound: if $n\ge9$, the squared distances of such a set are in ratio $(\kappa-1):\kappa$ for an integer $\kappa\ge3$ with $(2\kappa-1)^2\le n+1$ (Proposition~\ref{prop:tight}). In particular there is no such set for $9\le n\le23$, which gives $g(23)\le299$, one less than Blokhuis' bound (Corollary~\ref{cor:g23}). Finally, for every $d\ge5$ we describe all isosceles sets in $\R^d$ with at least $\binom{d+1}{2}+1$ points (Theorem~\ref{thm:large}).

The paper is organised as follows. Section~\ref{sec:decomp} contains the decomposition theory of isosceles sets: a short self-contained proof of Blokhuis' decomposition theorem, whose key step is a special case of Gallai's theorem on edge colourings of complete graphs without rainbow triangles, and the resulting general upper bound (Proposition~\ref{prop:bound}) in terms of $g$ and of the spherical analogue $M$. Section~\ref{sec:inputs} collects the facts about $g$ and $M$ that we use. Section~\ref{sec:proofs} proves Theorem~\ref{thm:main} for $10\le d\le22$ and Theorem~\ref{thm:structure} for $11\le d\le22$, Section~\ref{sec:nine} treats dimensions $9$ and $10$, and Section~\ref{sec:high} proves Theorem~\ref{thm:23}. Section~\ref{sec:formula} proves Theorem~\ref{thm:formula}. Section~\ref{sec:tight} deals with two-distance sets attaining Blokhuis' bound, and Section~\ref{sec:nopyr} proves Theorems~\ref{thm:nopyr} and~\ref{thm:all}. Section~\ref{sec:large} describes the large isosceles sets in every dimension, and Section~\ref{sec:open} lists the questions that remain open. The proof of Theorem~\ref{thm:all} in Sections~\ref{sec:formula}--\ref{sec:nopyr} uses from the earlier sections only Section~\ref{sec:decomp}; in particular it does not depend on \eqref{eq:g}.

\section{Decomposition of isosceles sets}\label{sec:decomp}

Throughout, $S$ is a finite set in $\R^d$, $\aff S$ is its affine hull and $\dim S=\dim\aff S$. For a distance $c$ occurring in $S$ let $\Gamma_c(S)$ be the graph on $S$ whose edges are the pairs at distance $c$. Colouring the pairs of $S$ by their distances, $S$ is isosceles if and only if no triangle receives three distinct colours. Edge colourings of complete graphs without such \emph{rainbow} triangles were studied by Gallai \cite{Gal67} (see \cite{GS04} for a modern account) and are now called Gallai colourings. The following lemma is a consequence of Gallai's structure theorem for such colourings, and it was also proved directly by Blokhuis \cite[Lemma~7.2.4]{Blo83}. We include a short proof to keep the paper self-contained.

\begin{lemma}\label{lem:gallai}
Let the edges of the complete graph on a finite vertex set $V$, $|V|\ge2$, be coloured so that no triangle receives three distinct colours. If at least three colours are used, then for some colour $c$ the spanning subgraph $\Gamma_c$ formed by the edges of colour $c$ is disconnected.
\end{lemma}
\begin{proof}
Induction on $n=|V|$. For $n\le3$ a colouring with three colours contains a rainbow triangle, so there is nothing to prove. Let $n\ge4$ and suppose that every used colour $c$ has $\Gamma_c$ connected. Then, as each $\Gamma_c$ is a connected spanning subgraph, every vertex is incident with an edge of every colour. Fix a vertex $v$ and let $V'=V\setminus\{v\}$. If the colouring restricted to $V'$ used at most two colours, some colour $b$ would occur only on edges at $v$; then $\Gamma_b$ would be a star centred at $v$, and being connected and spanning it would contain all edges at $v$, so that $v$ would be incident with edges of only one colour. Hence $V'$ carries at least three colours, and by induction there is a colour $c$ such that $\Gamma_c[V']$ is disconnected. Let $W_1,\dots,W_m$, $m\ge2$, be its components.

(i) Every $W_i$ contains a neighbour of $v$ in $\Gamma_c$: otherwise $W_i$ would be a union of components of $\Gamma_c$ not containing $v$, and $\Gamma_c$ would be disconnected.

(ii) For $i\ne j$, all edges between $W_i$ and $W_j$ have one and the same colour $e_{ij}$, and $e_{ij}\ne c$. Indeed, let $u\in W_i$, and let $x,x'\in W_j$ be joined by an edge of colour $c$. As $u$ and $x$ lie in different components of $\Gamma_c[V']$, neither $ux$ nor $ux'$ has colour $c$, and the triangle $uxx'$ forces $ux$ and $ux'$ to have the same colour. Since $W_j$ is connected in $\Gamma_c[V']$, the colour of $ux$ does not depend on $x\in W_j$ (this is trivial if $|W_j|=1$); by symmetry it does not depend on $u\in W_i$ either. Finally $e_{ij}\ne c$ because $W_i$ and $W_j$ are different components.

(iii) Let $z\in W_j$ be a neighbour of $v$ whose edge $vz$ has a colour $\delta\ne c$. For $i\ne j$ let $w_i\in W_i$ be a neighbour of $v$ in $\Gamma_c$, which exists by (i). The triangle $vzw_i$ has colours $\delta$, $c$ and $e_{ij}$, and $\delta\ne c\ne e_{ij}$; hence $e_{ij}=\delta$.

Now let $vz$ and $vz'$ be two edges at $v$ of colours $\delta,\delta'\ne c$, with $z\in W_j$ and $z'\in W_{j'}$. If $j=j'$, then (iii) gives $\delta=e_{ij}=\delta'$ for any $i\ne j$. If $j\ne j'$, then (iii) gives $\delta=e_{j'j}=e_{jj'}=\delta'$. So all edges at $v$ of colour different from $c$ have the same colour, and $v$ is incident with edges of at most two colours, a contradiction.
\end{proof}

\begin{lemma}[Blokhuis]\label{lem:decomp}
Let $S$ be a finite isosceles set in a metric space in which at least three distinct nonzero distances occur. Then there is a partition $S=X\sqcup Y$ with $|X|\ge2$ and $|Y|\ge1$ such that every $y\in Y$ is equidistant from all points of $X$.
\end{lemma}
\begin{proof}
Colour the pairs of $S$ by their distances. By Lemma~\ref{lem:gallai} there is a distance $c$ such that $\Gamma_c(S)$ is disconnected. Let $X$ be a component of $\Gamma_c(S)$ containing an edge; then $|X|\ge2$, and $Y=S\setminus X\ne\emptyset$. Let $y\in Y$ and let $xx'$ be an edge of $\Gamma_c(S)$ with $x,x'\in X$. Since $y$ is not in the component $X$, neither $d(y,x)$ nor $d(y,x')$ equals $c$, and the isosceles triangle $yxx'$ forces $d(y,x)=d(y,x')$. As $X$ is connected in $\Gamma_c(S)$, $d(y,\cdot)$ is constant on $X$.
\end{proof}

\begin{lemma}\label{lem:orth}
Let $X,Y\subset\R^d$ be finite nonempty sets such that every $y\in Y$ is equidistant from all points of $X$. Then $\ip{y-y'}{x-x'}=0$ for all $x,x'\in X$ and $y,y'\in Y$. Consequently $\dim(X\cup Y)\ge\dim X+\dim Y$, and if $|X|\ge2$ then $Y$ is contained in an affine subspace of dimension $d-\dim X$, namely the set of points equidistant from all points of $X$, which is orthogonal to $\aff X$ and meets it in the circumcentre $c_X$ of $X$; moreover $\nr{y-x}=\nr{y-c_X}+\rho_X^2$ for $x\in X$, $y\in Y$, where $\rho_X>0$ is the circumradius of $X$.
\end{lemma}
\begin{proof}
From $\nr{y-x}=\nr{y-x'}$ we get $2\ip{y}{x'-x}=\nr{x'}-\nr{x}$; subtracting the same identity for $y'$ gives $\ip{y-y'}{x-x'}=0$. Thus the direction spaces of $\aff X$ and $\aff Y$ are orthogonal, hence independent, and both lie in the direction space of $\aff(X\cup Y)$, which gives the dimension inequality. If $|X|\ge2$, the set $E$ of points of $\R^d$ equidistant from all of $X$ is the intersection of the bisecting hyperplanes of the pairs of points of $X$, an affine subspace whose direction space is the orthogonal complement of that of $\aff X$; it meets $\aff X$ in exactly one point $c_X$, which is equidistant from $X$, i.e.\ the circumcentre of $X$ in $\aff X$. For $y\in E$ and $x\in X$ the vectors $y-c_X$ and $c_X-x$ are orthogonal, whence the last formula.
\end{proof}

Following Ionin \cite[Definition~4.1]{Ion09}, a \emph{complete decomposition} of a finite set $S$ in a metric space is a partition $(S_1,\dots,S_k)$ of $S$ such that each $S_i$ is a two-distance set (at most two nonzero distances), $|S_i|\ge2$ for $i<k$, $|S_k|\ge1$, and for $i<j$ every point of $S_j$ is equidistant from all points of $S_i$.

\begin{proposition}[Ionin \cite{Ion09}]\label{prop:complete}
Every finite isosceles set $S$ in a metric space admits a complete decomposition. If $S\subset\R^d$ and $(S_1,\dots,S_k)$ is a complete decomposition, then $\dim S\ge\sum_i\dim S_i$.
\end{proposition}
\begin{proof}
Induction on $|S|$. If $S$ is a two-distance set, $(S)$ is a complete decomposition. Otherwise apply Lemma~\ref{lem:decomp} and choose $X$ of the smallest possible size. Then $X$ is a two-distance set: if $X$ had three distances, Lemma~\ref{lem:decomp} applied to $X$ would give $X=X_1\sqcup Y_1$, and $(X_1,\,Y_1\cup Y)$ would satisfy Lemma~\ref{lem:decomp} for $S$ with $|X_1|<|X|$ (points of $Y$ are equidistant from $X\supset X_1$). Since $|Y|<|S|$, $Y$ has a complete decomposition $(S_2,\dots,S_k)$ by induction, and $(X,S_2,\dots,S_k)$ is a complete decomposition of $S$. The dimension inequality follows from Lemma~\ref{lem:orth} and induction.
\end{proof}

For an integer $a\ge0$ let $g(a)$ be the largest size of a two-distance set in $\R^a$, with $g(0)=1$, and for $a\ge1$ let $M(a)$ be the largest size of a two-distance set in $\R^a$ that lies on a sphere; $M(1)=2$, and $M(a)\le g(a)$. Both functions are non-decreasing.

\begin{proposition}\label{prop:bound}
For every $d\ge1$,
\begin{equation}\label{eq:B}
\iso(d)\le\B(d):=\max\left\{\sum_{i=1}^{k-1}M(a_i)+g(a_k)\ :\ \begin{array}{l}k\ge1,\ a_1,\dots,a_{k-1}\ge1,\\[2pt] a_k\ge0,\ \sum_{i=1}^k a_i\le d\end{array}\right\}.
\end{equation}
\end{proposition}
\begin{proof}
Let $S\subset\R^d$ be isosceles with complete decomposition $(S_1,\dots,S_k)$ and put $a_i=\dim S_i$. Then $\sum a_i\le\dim S\le d$ by Proposition~\ref{prop:complete}, $a_i\ge1$ for $i<k$ since $|S_i|\ge2$, and $|S_k|\le g(a_k)$. For $i<k$ the set $S_i$ lies on a sphere of $\R^d$ centred at any point of $S_k\ne\emptyset$; the intersection of this sphere with $\aff S_i$ is a sphere of $\aff S_i$ (of positive radius, as $|S_i|\ge2$), so $|S_i|\le M(a_i)$.
\end{proof}

\begin{remark}\label{rem:blokhuis}
Blokhuis \cite[Theorem~7.2.5]{Blo83} obtained from the decomposition the bound $\iso(d)\le\binom{d+2}{2}$, with equality only for two-distance sets and for spherical two-distance sets together with their centre, by using in \eqref{eq:B} the general bounds $g(a)\le\binom{a+2}{2}$ \cite{Blo84} and $M(a)\le a(a+3)/2$ \cite{DGS77}. Proposition~\ref{prop:bound} is the same argument with the exact values of $g$ and $M$ where they are known.
\end{remark}

\begin{remark}\label{rem:lower}
If $X$ is a spherical two-distance set with centre $o$ and radius $r$, then $X\cup\{o\}$ is isosceles, and so is $X\cup\{o,o+rv,o-rv\}$ for a unit vector $v$ orthogonal to $\aff X$ (in the triangle $o,o+rv,x$ two sides have length $r$; in all other triangles two sides are equal because $o$, $o\pm rv$ are equidistant from $X$ and $o$ is equidistant from $o\pm rv$). Hence
\begin{equation}\label{eq:lower}
\iso(d)\ge\max\{g(d),\ M(d)+1,\ M(d-1)+3\}.
\end{equation}
Kelly's $8$-point set in $\R^3$ and Ionin's $30$-point set in $\R^7$ are of the third type. The lower bound $\binom{d+1}{2}+1$ is also recorded on the problem page \cite{EP503}.
\end{remark}

The next lemma describes a two-distance set $T$ all of whose points are equidistant from a set $X$, when $T$ together with a point of $X$ is isosceles. It is used in Sections~\ref{sec:nine} and~\ref{sec:formula}.

\begin{lemma}\label{lem:levels}
Let $X,T\subset\R^d$ be finite sets with $|X|\ge2$ and $T\ne\emptyset$ such that every point of $T$ is equidistant from all points of $X$, and put $m=c_X$ and $\rho=\rho_X$. Suppose that $T$ is a two-distance set, let $\Delta$ be the set of squared distances between distinct points of $T$, and suppose that $T\cup\{p\}$ is isosceles for some $p\in X$. Let $C$ be the set of all $x\in T$ with $\rho^2+\nr{x-m}\notin\Delta$. Then
\begin{enumerate}[label=\textup{(\roman*)}]
\item all points of $C$ have the same distance from $m$;
\item every $y\in T\setminus C$ satisfies $\nr{y-x}=\nr{y-p}$ for all $x\in C$, so that $y$ is equidistant from all points of $C\cup X$;
\item if $C=\emptyset$, then $T\cup\{p\}$ is a two-distance set.
\end{enumerate}
\end{lemma}
\begin{proof}
By Lemma~\ref{lem:orth},
\[
\nr{p-x}=\rho^2+r(x)\qquad(x\in T),\qquad r(x):=\nr{x-m}.
\]
For distinct $x,y\in T$ the triangle $p,x,y$ is isosceles, so
\begin{equation}\label{eq:iso}
r(x)=r(y)\quad\text{or}\quad\nr{x-y}\in\{\rho^2+r(x),\ \rho^2+r(y)\}.
\end{equation}
Let $x,y\in C$ be distinct. Then $\nr{x-y}\in\Delta$, while $\rho^2+r(x)$ and $\rho^2+r(y)$ are not in $\Delta$, so \eqref{eq:iso} gives $r(x)=r(y)$. This proves (i). Let $x\in C$ and $y\in T\setminus C$. Then $\rho^2+r(y)\in\Delta$ and $\rho^2+r(x)\notin\Delta$, so $r(x)\ne r(y)$ and $\nr{x-y}\ne\rho^2+r(x)$, and \eqref{eq:iso} gives $\nr{x-y}=\rho^2+r(y)=\nr{p-y}$. This proves (ii), because $y$ is equidistant from all points of $X$. If $C=\emptyset$, then $\nr{p-x}\in\Delta$ for all $x\in T$, so all squared distances in $T\cup\{p\}$ lie in $\Delta$, which has at most two elements. This proves (iii).
\end{proof}

\section{Two-distance sets: the inputs}\label{sec:inputs}

The values of $g$ in \eqref{eq:g} and the following facts are all that we need from the theory of two-distance sets. For a set whose squared distances lie in $\{2,4\}$, the \emph{short-distance graph} is the graph on the set whose edges are the pairs at squared distance $2$.

\begin{theorem}[\cite{ES66,Lis97,BLSTYY}]\label{thm:g}
Let $1\le d\le22$. Every two-distance set in $\R^d$ has at most $g(d)$ points, with $g$ as in \eqref{eq:g}. For $10\le d\le21$ every two-distance set with $g(d)$ points in $\R^d$ is similar to $\J(d+1,2)$. Every two-distance set with $45$ points and affine dimension $9$ is similar to $\J(10,2)$.
\end{theorem}

For $d\ge9$ this is \cite[Theorem~1.1 and Corollary~1.2]{BLSTYY}. For the last statement recall the theorem of Larman, Rogers and Seidel \cite{LRS77}: if a two-distance set in $\R^d$ has more than $2d+3$ points, then the ratio of its squared distances is $(\kappa-1)/\kappa$ for an integer $\kappa$ with $2\le\kappa\le(1+\sqrt{2d})/2$. For $d\le12$ this forces $\kappa=2$, so a $45$-point two-distance set of affine dimension $9$ has squared distances $2$ and $4$ after scaling, and \cite[Theorem~1.1]{BLSTYY} shows that it is similar to $\J(10,2)$.

\begin{theorem}[\cite{BLSTYY}]\label{thm:conn}
Let $X$ be a finite set with squared distances in $\{2,4\}$, affine dimension $d\ge9$ and connected short-distance graph. If $|X|\ge\binom{d+1}{2}-1$, then $X$ is congruent to a subset of $\J(d+1,2)$.
\end{theorem}

This is contained in the proof of \cite[Proposition~2.1]{BLSTYY}. There a point of $X$ is translated to the origin, and $X$ then lies in an even lattice spanned by vectors of squared norm $2$. Its root system is irreducible of type $A_d$ or $D_d$, or it is reducible (an irreducible root system of type $E$ has rank at most $8$). In the reducible case and in type $D_d$ the proof gives $|X|\le\binom{d+1}{2}-2$. In type $A_d$ it shows that, after a translation and possibly the map $x\mapsto\mathbf 1-x$, the set $X$ is a family of $k$-subsets of $\{1,\dots,d+1\}$, identified with their characteristic vectors, and that $|X|\le\binom{d+1}{2}-2$ unless $k=2$. For $k=2$ the set $X$ is a subset of $\J(d+1,2)$.

\begin{theorem}[\cite{Mus09,DGS77}]\label{thm:M}
$M(d)=d(d+1)/2$ for $7\le d\le21$, and in particular $M(7)=28$ and $M(8)=36$. Moreover $M(22)=275$, and up to isometry the McLaughlin configuration is the only spherical two-distance set with $275$ points in $\R^{22}$.
\end{theorem}

Musin proves $M(d)=d(d+1)/2$ for $6<d<22$ and $23<d<40$ \cite[Theorem~7]{Mus09}. Barg and Yu \cite{BY13} settled $d=23$, with $M(23)=276$, and Glazyrin and Yu \cite{GY18} proved $M(d)=d(d+1)/2$ for all $d\ge7$ except possibly $d=(2k+1)^2-3$. For $d=22$ the absolute bound $d(d+3)/2$ of Delsarte, Goethals and Seidel \cite{DGS77} gives $M(22)\le275$, the McLaughlin configuration attains it \cite[Example~4.10]{DGS77}, and the extremal set is unique up to isometry \cite[Remark~4]{Mus09}. For $d\le6$ we only use $M(d)\le g(d)$, which is attained for $2\le d\le6$: $M(2)=5$ (pentagon), $M(3)=6$ (octahedron), $M(4)=10$, $M(5)=16$, $M(6)=27$ (the sets $\J(5,2)$, the $16$-point half-cube and the $27$-point Schl\"afli configuration are spherical \cite{Lis97}). Table~\ref{tab:values} lists these numbers.

\begin{table}[ht]
\centering\small
\begin{tabular}{l*{12}{r}}
\toprule
$d$ & 1&2&3&4&5&6&7&8&9&10&11&12\\
\midrule
$g(d)$ & 3&5&6&10&16&27&29&45&45&55&66&78\\
$M(d)$ & 2&5&6&10&16&27&28&36&45&55&66&78\\
$\B(d)$ & 3&6&8&11&17&28&30&45&47&56&67&79\\
$\iso(d)$ & 3&6&8&11&17&28&30&45&46&56&67&79\\
\midrule
$d$ & 13&14&15&16&17&18&19&20&21&22&23&\\
\midrule
$g(d)$ & 91&105&120&136&153&171&190&210&231&275&$277$--$299$&\\
$M(d)$ & 91&105&120&136&153&171&190&210&231&275&276&\\
$\B(d)$ & 92&106&121&137&154&172&191&211&232&276&$\le299$&\\
$\iso(d)$ & 92&106&121&137&154&172&191&211&232&276&$278$--$299$&\\
\bottomrule
\end{tabular}
\medskip
\caption{The largest two-distance sets ($g$), spherical two-distance sets ($M$), the bound $\B$ of \eqref{eq:B}, and the largest isosceles sets ($\iso$) in $\R^d$. The values $\iso(d)$ for $d\le8$ are due to Kelly, Croft, Kido and Ionin. The values for $9\le d\le22$ and the bounds for $d=23$ are proved here. Entries of the form $a$--$b$ are the known lower and upper bounds, and for $d=23$ the bound $\B(23)\le299$ uses $g(23)\le299$ (Corollary~\ref{cor:g23}).}\label{tab:values}
\end{table}

\section{Dimensions ten to twenty-two}\label{sec:proofs}

\begin{lemma}\label{lem:B}
For $10\le d\le21$ the maximum in \eqref{eq:B} equals $\binom{d+1}{2}+1$, and it is attained only for $k=2$, $(a_1,a_2)=(d,0)$. The second largest value of the expression in \eqref{eq:B} is $\binom{d+1}{2}$, attained only for $k=1$, $a_1=d$. For $d=22$ the same holds with the values $276$ and $275$. For $11\le d\le22$ the third largest value is $\binom d2+3$, attained only for the compositions $(d-1,1)$, $(d-1,1,0)$ and $(1,d-1,0)$. For $d=10$ the values following $55$ are $50$, from $(2,8)$, then $49$, from $(1,1,8)$, and then $48=\binom{10}{2}+3$, from $(9,1)$, $(9,1,0)$ and $(1,9,0)$. For $d=9$ the three largest values are $47$, attained only for $(1,8)$, then $46$, attained only for $(9,0)$, and then $45$, from $(9)$ and $(8)$.
\end{lemma}
\begin{proof}
This is a finite computation with the values of Table~\ref{tab:values}, which we carried out by enumerating all compositions $(a_1,\dots,a_k)$ with $\sum a_i\le d$. It can also be checked by hand: for $a\le21$ one has $M(a)\le\binom{a+1}{2}+\epsilon(a)$ and $g(a)\le\binom{a+1}{2}+\epsilon(a)$ with $\epsilon(a)=0$ except $\epsilon(0)=1$, $\epsilon(1)=2$, $\epsilon(2)=2$, $\epsilon(5)=1$, $\epsilon(6)=6$, $\epsilon(7)=1$, $\epsilon(8)=9$, and $\binom{x+1}{2}$ is superadditive, so that for $k\ge2$ the sum in \eqref{eq:B} is at most $\binom{d+1}{2}-\sum_{i<j}a_ia_j+\sum_i\epsilon(a_i)$. A composition with value at least $\binom d2+3=\binom{d+1}{2}-d+3$ therefore satisfies $\sum_{i<j}a_ia_j\le\sum_i\epsilon(a_i)+d-3$, and there are few of these. For $d=22$ the only compositions with a part equal to $22$ are $(22)$ and $(22,0)$, and all other compositions have parts at most $21$, so the same estimate applies to them.
\end{proof}

\begin{proof}[Proof of Theorem~\ref{thm:main} for $10\le d\le22$]
The lower bound follows from \eqref{eq:lower} and the values $M(d)=\binom{d+1}{2}$ for $d\le21$ and $M(22)=275$. Let $S\subset\R^d$ be isosceles with $|S|\ge\B(d)$, and let $(S_1,\dots,S_k)$ be a complete decomposition. By the proof of Proposition~\ref{prop:bound} and Lemma~\ref{lem:B}, $k=2$, $\dim S_1=d$, $S_2=\{o\}$ is a single point, and $|S_1|=M(d)$. Since $o$ is equidistant from $S_1$, the set $S_1$ lies on a sphere centred at $o$, and $o$ is the centre of the circumsphere of $S_1$, which is unique because $S_1$ spans $\R^d$. For $d\le21$, $S_1$ is a two-distance set with $g(d)$ points, hence similar to $\J(d+1,2)$ by Theorem~\ref{thm:g}, and the circumcentre of $\J(d+1,2)$ is its centroid. For $d=22$, $S_1$ is a spherical two-distance set with $275$ points, hence the McLaughlin configuration by Theorem~\ref{thm:M}.
\end{proof}

\begin{proof}[Proof of Theorem~\ref{thm:structure} for $11\le d\le22$]
Let $S$ be isosceles with more than $\binom d2+3$ points and $(S_1,\dots,S_k)$ a complete decomposition. By Lemma~\ref{lem:B} either $k=1$, so that $S$ is a two-distance set, or $k=2$ with $a_2=0$, so that $S=S_1\cup\{o\}$ with $S_1$ a two-distance set lying on a sphere centred at $o$. The example is Remark~\ref{rem:lower} with $X=\J(d,2)$, which has $\binom d2$ points, spans a hyperplane of $\R^d$ and has squared distances $2$, $4$ and squared circumradius $\rho^2=2-4/d$. The additional distances $\rho$, $\sqrt2\rho$ and $2\rho$ are different from $\sqrt2$ and $2$, so the set has at least three distances, and removing any one of its points leaves a set with at least three distances, so it is of neither type. (It has the complete decompositions $(\J(d,2),\{o,o\pm rv\})$, $(\J(d,2),\{o\pm rv\},\{o\})$ and $(\{o\pm rv\},\J(d,2),\{o\})$, corresponding to the compositions $(d-1,1)$, $(d-1,1,0)$ and $(1,d-1,0)$ of Lemma~\ref{lem:B}.)
\end{proof}

\section{Dimension nine}\label{sec:nine}

In $\R^9$ the bound of Proposition~\ref{prop:bound} is $47$, coming from the composition $(1,8)$: two points $S_1=\{p,q\}$ and a $45$-point two-distance set $S_2$ of dimension $8$, every point of which is equidistant from $p$ and $q$. Lison\v{e}k's set $\LL$ is a $45$-point two-distance set of dimension $8$, so this composition cannot be dismissed by counting. Lemma~\ref{lem:noext} below shows that no isosceles set contains such a configuration, in any dimension, and that the same holds when the last part has $44$ points. Together with Lemma~\ref{lem:B} this proves Theorem~\ref{thm:main} for $d=9$, and it is also needed for Theorem~\ref{thm:structure} in dimensions $9$ and $10$, where the compositions $(1,8)$, $(2,8)$ and $(1,1,8)$, ending with an $8$-dimensional part, have values above the thresholds $45$ and $48$. We first record some properties shared by all large two-distance sets of dimension $8$.

\begin{lemma}\label{lem:eight}
Let $T$ be a two-distance set of affine dimension $8$ with $|T|\ge44$. Then the ratio of the squared distances of $T$ is $2$. Scale $T$ so that its squared distances are $2$ and $4$, and for $T_0\subseteq T$ let $\Gamma(T_0)$ be the short-distance graph of $T_0$.
\begin{enumerate}[label=\textup{(\roman*)}]
\item Every $T_0\subseteq T$ with $|T_0|\ge43$ has affine dimension $8$.
\item No $T_0\subseteq T$ with $|T_0|\ge37$ lies on a sphere.
\item If $T_0\subseteq T$ has a partition $T_0=C\sqcup D$ with $|C|\ge2$ and $D\ne\emptyset$ such that every point of $D$ is equidistant from all points of $C$, then $|T_0|\le37$.
\item $\Gamma(T_0)$ is connected for every $T_0\subseteq T$ with $|T_0|\ge38$.
\item $\Gamma(T)$ is prime: every set $C\subseteq T$ such that each point of $T\setminus C$ is adjacent in $\Gamma(T)$ to all points of $C$ or to none of them satisfies $|C|\le1$ or $C=T$.
\end{enumerate}
\end{lemma}
\begin{proof}
Since $|T|>2\cdot8+3$, the Larman--Rogers--Seidel theorem \cite{LRS77} gives the ratio $(\kappa-1)/\kappa$ with $2\le\kappa\le(1+\sqrt{16})/2$, so $\kappa=2$.

(i) A two-distance set of affine dimension at most $7$ has at most $g(7)=29$ points.

(ii) A spherical two-distance set of affine dimension $a\le8$ has at most $M(a)\le M(8)=36$ points.

(iii) Let $a=\dim C$ and $b=\dim D$. By Lemma~\ref{lem:orth}, $a+b\le\dim T_0\le8$, and $a\ge1$ because $|C|\ge2$. The set $C$ lies on a sphere centred at any point of $D$, and hence on a sphere of positive radius in $\aff C$. Therefore $|T_0|\le M(a)+g(b)$. By Table~\ref{tab:values} the largest value of $M(a)+g(b)$ with $a\ge1$ and $a+b\le8$ is $M(8)+g(0)=37$, the next ones being $M(2)+g(6)=M(6)+g(2)=32$.

(iv) Suppose that $\Gamma(T_0)$ is disconnected. If all its components are single points, then all distances in $T_0$ are equal and $|T_0|\le9$. Otherwise let $C$ be a component with at least two points and $D=T_0\setminus C\ne\emptyset$. Every point of $D$ is at squared distance $4$ from all points of $C$, and (iii) gives $|T_0|\le37$.

(v) A set $C$ as in (v) with $2\le|C|<|T|$ satisfies the hypothesis of (iii) with $T_0=T$ and $D=T\setminus C$, because each point of $D$ is at squared distance $2$ from all points of $C$ or at squared distance $4$ from all of them. So $|T|\le37$, a contradiction.
\end{proof}

\begin{lemma}\label{lem:noext}
Let $S\subset\R^d$ be an isosceles set with a complete decomposition $(S_1,\dots,S_k)$, $k\ge2$. Then $S_k$ is not a two-distance set of affine dimension $8$ with at least $44$ points.
\end{lemma}
\begin{proof}
Suppose that $T=S_k$ is such a set, scaled so that its squared distances are $2$ and $4$ (Lemma~\ref{lem:eight}). Let $X=S_1$, so $|X|\ge2$ and every point of $T$ is equidistant from $X$, and fix $p\in X$. We apply Lemma~\ref{lem:levels} with $\Delta=\{2,4\}$, $m=c_X$ and $\rho=\rho_X$. Let $C$ be the set of all $x\in T$ with $\rho^2+\nr{x-m}\notin\{2,4\}$, so that
\[
\nr{p-x}\in\{2,4\}\qquad(x\in T\setminus C)
\]
by Lemma~\ref{lem:orth}. By Lemma~\ref{lem:levels}(ii) every point of $T\setminus C$ is equidistant from all points of $C$, and by Lemma~\ref{lem:eight}(iii), $|C|\le1$ or $C=T$. If $C=T$ then $T$ lies on a sphere about $m$ by Lemma~\ref{lem:levels}(i), contradicting Lemma~\ref{lem:eight}(ii).

Note that $p\notin\aff T$: by Lemma~\ref{lem:orth}, $T$ lies in the affine subspace of points equidistant from $X$, which meets $\aff X$ only in $m$, and $p\ne m$. Hence $\dim(T'\cup\{p\})=\dim T'+1=9$ for every $T'\subseteq T$ with $\aff T'=\aff T$. We also use that $\J(10,2)$ minus one or two points has affine dimension $9$: for distinct $i,j$, at least six of the eight differences $(e_i+e_k)-(e_j+e_k)$ survive, so the remaining differences span $\{x\in\R^{10}:\sum_ix_i=0\}$.

\emph{Case $C=\emptyset$.} Then $\nr{p-x}\in\{2,4\}$ for all $x\in T$, so $T\cup\{p\}$ is a two-distance set with $|T|+1\ge45$ points and affine dimension $9$. Since $g(9)=45$, this forces $|T|=44$, and by Theorem~\ref{thm:g} the set $T\cup\{p\}$ is similar to $\J(10,2)$. Then $T$ is similar to $\J(10,2)$ minus a point, which has affine dimension $9$, whereas $\dim T=8$.

\emph{Case $|C|=1$, $C=\{c\}$.} Then $T'=(T\setminus\{c\})\cup\{p\}$ is a two-distance set with $|T|\ge44$ points, and it has affine dimension $9$ because $\aff(T\setminus\{c\})=\aff T$ by Lemma~\ref{lem:eight}(i). If $p$ has no neighbour in the short-distance graph of $T'$, then $T\setminus\{c\}$ lies on the sphere of squared radius $4$ about $p$, contradicting Lemma~\ref{lem:eight}(ii). Otherwise the short-distance graph of $T'$ is connected, because $\Gamma(T\setminus\{c\})$ is connected by Lemma~\ref{lem:eight}(iv). As $|T'|\ge44=\binom{10}{2}-1$, Theorem~\ref{thm:conn} shows that $T'$ is congruent to a subset of $\J(10,2)$. Then $T\setminus\{c\}=T'\setminus\{p\}$ is congruent to $\J(10,2)$ minus one or two points, which has affine dimension $9$, contradicting $\dim(T\setminus\{c\})=8$.
\end{proof}

\begin{proof}[Proof of Theorems~\ref{thm:main} and \ref{thm:structure} for $d=9$]
The lower bound $\iso(9)\ge46$ is \eqref{eq:lower}. Let $S\subset\R^9$ be isosceles with $|S|\ge46$ and let $(S_1,\dots,S_k)$ be a complete decomposition with dimensions $a_i$. By Lemma~\ref{lem:B} the composition is $(1,8)$ or $(9,0)$. In the first case $|S_1|\le M(1)=2$, so $S_2$ is a two-distance set of dimension $8$ with at least $44$ points, which Lemma~\ref{lem:noext} excludes. In the second case $|S_1|=45$, and $S_1$ spans $\R^9$ and lies on a sphere centred at the point of $S_2$. By Theorem~\ref{thm:g}, $S_1$ is similar to $\J(10,2)$, and $S_2$ is its centroid. This proves both theorems for $d=9$.
\end{proof}

\begin{proof}[Completion of the proof of Theorem~\ref{thm:structure} for $d=10$]
By Lemma~\ref{lem:B} the values of \eqref{eq:B} for $d=10$ are $56$, $55$, $50$, $49$, $48,\dots$. The values $50$ and $49$ come from the compositions $(2,8)$ and $(1,1,8)$, whose last part is an $8$-dimensional two-distance set. If the last part has at least $44$ points, Lemma~\ref{lem:noext} excludes the decomposition, and if it has at most $43$ points these compositions give at most $5+43=48$. Hence an isosceles set in $\R^{10}$ with more than $48=\binom{10}{2}+3$ points has a decomposition of type $(10)$ or $(10,0)$, as claimed, and the example of Section~\ref{sec:proofs} shows that $48$ is attained.
\end{proof}

\section{Dimension 23}\label{sec:high}

\begin{proof}[Proof of Theorem~\ref{thm:23}]
We have $277\le g(23)\le299$ by \cite{GKM26} and by Corollary~\ref{cor:g23}, whose proof does not depend on the present section. (Blokhuis' bound gives $g(23)\le\binom{25}{2}=300$. By the Larman--Rogers--Seidel theorem \cite{LRS77}, the squared distances of a two-distance set with more than $2\cdot23+3$ points in $\R^{23}$ are in ratio $(\kappa-1):\kappa$ with $\kappa\in\{2,3\}$, and if $\kappa=2$ the set has at most $\binom{24}{2}=276$ points by \cite[Theorem~1.1]{BLSTYY}; so a two-distance set with more than $276$ points in $\R^{23}$ has squared distances in ratio $2:3$. For this ratio the relative bound \cite[Proposition~4.1]{BLSTYY} and the bound of Chen and Yu \cite[Theorem~3.1]{CY25} give $575$ and $577$, which are weaker than Blokhuis' bound.)

Now evaluate \eqref{eq:B} for $d=23$. The value $g(23)$ occurs only for the composition $(23)$. Every other composition involves only $M(a)$ with $a\le23$ and $g(a)$ with $a\le22$, whose values are given in Table~\ref{tab:values}, and the enumeration of Lemma~\ref{lem:B} shows that these compositions have value at most $278=M(22)+g(1)$, attained by $(22,1)$, $(22,1,0)$ and $(1,22,0)$. Hence $\iso(23)\le\max\{g(23),278\}$. Conversely $\iso(23)\ge g(23)$, and \eqref{eq:lower} with the McLaughlin configuration gives $\iso(23)\ge M(22)+3=278$. The bounds $278\le\iso(23)\le299$ follow from the first paragraph.
\end{proof}

\begin{remark}\label{rem:forum}
In every dimension where $\iso(d)$ is known, including $9\le d\le22$ by Theorem~\ref{thm:main}, the lower bound \eqref{eq:lower} is the truth (Table~\ref{tab:values}), in agreement with the identity \eqref{eq:forum} stated in \cite{EP503} for all $d\ge2$. In dimension $23$ the identity \eqref{eq:forum} reads $\iso(23)=\max\{g(23),278\}$, which is Theorem~\ref{thm:23}. Dimension $9$ shows that the bound \eqref{eq:B} alone does not always give \eqref{eq:forum}: there $\B(9)=47>46=\iso(9)$. In Sections~\ref{sec:formula}--\ref{sec:nopyr} we show that the composition $(1,d-1)$, which is responsible for this, is the only obstacle to \eqref{eq:forum} in any dimension, and that it never prevents \eqref{eq:forum}.
\end{remark}

\section{The general formula}\label{sec:formula}

In this section we prove Theorem~\ref{thm:formula}. We use the general bounds
\begin{equation}\label{eq:general}
g(a)\le\binom{a+2}{2},\qquad\binom{a+1}{2}\le M(a)\le\frac{a(a+3)}{2}
\end{equation}
of Blokhuis \cite{Blo84} and of Delsarte, Goethals and Seidel \cite{DGS77}, which were already used in Remark~\ref{rem:blokhuis}; the lower bound for $M(a)$ is given by the set $\J(a+1,2)$. Let $L(d)$ be the right-hand side of \eqref{eq:forum}. Then $\iso(d)\ge L(d)$ by \eqref{eq:lower}, and by \eqref{eq:general}
\begin{equation}\label{eq:Lbig}
L(d)\ \ge\ M(d)+1\ \ge\ \binom{d+1}{2}+1 .
\end{equation}

We call the sequence $(a_1,\dots,a_k)$ of the dimensions $a_i=\dim S_i$ the \emph{type} of a complete decomposition $(S_1,\dots,S_k)$. By the proof of Proposition~\ref{prop:bound}, if $S\subset\R^d$ has a complete decomposition of type $a=(a_1,\dots,a_k)$, then $a$ is a composition as in \eqref{eq:B}, that is, $a_i\ge1$ for $i<k$, $a_k\ge0$ and $\sum_ia_i\le d$, and
\begin{equation}\label{eq:parts}
|S_i|\le M(a_i)\quad(i<k),\qquad |S_k|\le g(a_k).
\end{equation}
For a composition $a=(a_1,\dots,a_k)$ put $|a|=\sum_ia_i$, $e_2(a)=\sum_{i<j}a_ia_j$ and
\[
v(a):=\sum_{i<k}\frac{a_i(a_i+3)}{2}+\binom{a_k+2}{2},
\]
so that $|S|\le v(a)$ by \eqref{eq:parts} and \eqref{eq:general}.

\begin{lemma}\label{lem:v}
$v(a)=\binom{|a|+2}{2}-e_2(a)$ for every composition $a$.
\end{lemma}
\begin{proof}
We have $\frac{x(x+3)}{2}=\binom{x+1}{2}+x$ and $\binom{x+2}{2}=\binom{x+1}{2}+x+1$, so $v(a)=\sum_i\binom{a_i+1}{2}+|a|+1$. Since $\binom{x+y+1}{2}=\binom{x+1}{2}+\binom{y+1}{2}+xy$, induction on $k$ gives $\sum_i\binom{a_i+1}{2}=\binom{|a|+1}{2}-e_2(a)$, and $\binom{|a|+1}{2}+|a|+1=\binom{|a|+2}{2}$.
\end{proof}

For $d\ge2$ we call the compositions
\[
(d),\quad(d,0),\quad(d-1,1),\quad(d-1,1,0),\quad(1,d-1,0),\quad(1,d-1)
\]
the \emph{special} compositions of $d$.

\begin{lemma}\label{lem:nonspecial}
Let $d\ge2$ and let $a$ be a composition with $|a|\le d$ which is not special. Then $v(a)\le\binom{d+1}{2}+1$, and $v(a)\le\binom{d+1}{2}$ if $d\ge5$.
\end{lemma}
\begin{proof}
If $|a|\le d-1$, then $v(a)\le\binom{d+1}{2}$ by Lemma~\ref{lem:v}. Let $|a|=d$, so that $v(a)=\binom{d+1}{2}+d+1-e_2(a)$. We show that $e_2(a)\ge d$, and that $e_2(a)\ge d+1$ if $d\ge5$. Let $b_1\ge\dots\ge b_r\ge1$ be the positive parts of $a$; then $\sum_jb_j=d$ and $e_2(a)=e_2(b)$. If $r=1$, then $a$ is $(d)$ or $(d,0)$. If $r=2$ and $b_2=1$, then $a$ is one of $(d-1,1)$, $(1,d-1)$, $(d-1,1,0)$, $(1,d-1,0)$. If $r=2$ and $b_2\ge2$, then $d\ge4$ and $e_2(a)=b_2(d-b_2)\ge2(d-2)$, which is at least $d$, and at least $d+1$ if $d\ge5$. If $r\ge3$, then $\sum_jb_j^2\le(d-r+1)^2+(r-1)$, hence
\[
e_2(a)=\tfrac12\Bigl(d^2-\sum_jb_j^2\Bigr)\ \ge\ \tfrac12(r-1)(2d-r)\ \ge\ 2d-3,
\]
because $\frac12(r-1)(2d-r)$ is non-decreasing in $r$ for $3\le r\le d$. Finally $2d-3\ge d$ as $d\ge r\ge3$, and $2d-3\ge d+1$ if $d\ge4$.
\end{proof}

\begin{lemma}\label{lem:types}
Let $d\ge2$ and let $S\subset\R^d$ be an isosceles set with a complete decomposition of type $a\ne(1,d-1)$. Then $|S|\le L(d)$. If $d\ge5$ and $a$ is not special, then $|S|\le\binom{d+1}{2}$.
\end{lemma}
\begin{proof}
If $a$ is not special, then $|S|\le v(a)$, and the claims follow from Lemma~\ref{lem:nonspecial} and \eqref{eq:Lbig}. If $a$ is special, then \eqref{eq:parts} with $g(0)=1$, $g(1)=3$ and $M(1)=2$ gives $|S|\le g(d)$ for $a=(d)$, $|S|\le M(d)+1$ for $a=(d,0)$, $|S|\le M(d-1)+g(1)=M(d-1)+3$ for $a=(d-1,1)$, and $|S|\le M(d-1)+M(1)+g(0)=M(d-1)+3$ for $a=(d-1,1,0)$ and $a=(1,d-1,0)$.
\end{proof}

So the only type that is not controlled by counting is $(1,d-1)$: two points together with a two-distance set of dimension $d-1$ in their bisecting hyperplane. For this type \eqref{eq:parts} gives only $|S|\le2+g(d-1)$, which exceeds $L(d)$ for $d=9$.

\begin{lemma}\label{lem:twopoints}
Let $S\subset\R^d$ be an isosceles set with a complete decomposition $(\{p,q\},T)$, and let $m$ be the midpoint of $p$ and $q$. Then at least one of the following holds.
\begin{enumerate}[label=\textup{(\alph*)}]
\item All points of $T$ have the same distance from $m$.
\item $S$ has a complete decomposition $(\{p,q\},C,T\setminus C)$ with three parts.
\item $S$ has a complete decomposition $(\{p,q,c\},T\setminus\{c\})$ whose first part has dimension $2$.
\item $T\cup\{p\}$ is a two-distance set, and hence a two-distance pyramid with apex $p$ and base $T$.
\end{enumerate}
\end{lemma}
\begin{proof}
The set $T$ is a nonempty two-distance set, all of its points are equidistant from $p$ and $q$, and $T\cup\{p\}$ is isosceles. So Lemma~\ref{lem:levels} applies to $X=\{p,q\}$, $T$ and $p$, with $c_X=m$; let $C\subseteq T$ be the set defined there.

If $C=T$, then (a) holds by Lemma~\ref{lem:levels}(i). If $C=\emptyset$, then $T\cup\{p\}$ is a two-distance set by Lemma~\ref{lem:levels}(iii). Its base $T$ lies in the bisecting hyperplane $H$ of $p$ and $q$, and $p\notin H$. This is (d).

Let $\emptyset\ne C\ne T$. By Lemma~\ref{lem:levels}(ii) every point of $T\setminus C$ is equidistant from $p$, $q$ and all points of $C$. If $|C|\ge2$, then $(\{p,q\},C,T\setminus C)$ is a complete decomposition of $S$: its parts are two-distance sets, the first two have at least two points, the third is nonempty, and all points of $C$ and of $T\setminus C$ are equidistant from $p$ and $q$. This is (b). If $C=\{c\}$, then $c\ne m$, since otherwise the three collinear points $p,m,q$ would lie on a sphere centred at a point of $T\setminus\{c\}$. So $p,q,c$ are not collinear, $\{p,q,c\}$ is a two-distance set because $\nr{p-c}=\nr{q-c}$, and $(\{p,q,c\},T\setminus\{c\})$ is a complete decomposition of $S$. This is (c).
\end{proof}

\begin{lemma}\label{lem:double}
Let $X\subset\R^d$ be a two-distance pyramid with apex $p$ and base $T\subset H$, and let $q$ be the mirror image of $p$ in the hyperplane $H$. Then $X\cup\{q\}$ is an isosceles set with $|X|+1$ points. Hence $\iso(d)\ge\pyr(d)+1$.
\end{lemma}
\begin{proof}
We have $q\notin X$ because $q\notin H$ and $q\ne p$. A triangle with vertices in $X$, or in $T\cup\{q\}$, which is congruent to $X$, has its three sides in a set of at most two lengths, so two of them are equal. The remaining triangles are $p,q,x$ with $x\in T$, and $\nr{p-x}=\nr{q-x}$.
\end{proof}

\begin{proof}[Proof of Theorem~\ref{thm:formula}]
The lower bound is \eqref{eq:lower} together with Lemma~\ref{lem:double}. For the upper bound let $S\subset\R^d$ be isosceles and let $(S_1,\dots,S_k)$ be a complete decomposition of $S$, of type $a$. If $a\ne(1,d-1)$, then $|S|\le L(d)$ by Lemma~\ref{lem:types}. Let $a=(1,d-1)$. Then $S_1$ consists of collinear points on a sphere centred at a point of $S_2$, so $S_1=\{p,q\}$ has two points, and $T:=S_2$ has dimension $d-1\ge1$. We apply Lemma~\ref{lem:twopoints}.

In case (a), $T$ lies on a sphere centred at the midpoint of $p$ and $q$, of positive radius because $|T|\ge2$; thus $|T|\le M(d-1)$ and $|S|\le M(d-1)+2<L(d)$. In cases (b) and (c), $S$ has a complete decomposition with three parts, or with a first part of dimension $2$; its type is not $(1,d-1)$, and $|S|\le L(d)$ by Lemma~\ref{lem:types}. In case (d), $T\cup\{p\}$ is a two-distance pyramid in $\R^d$, so $|S|=|T\cup\{p\}|+1\le\pyr(d)+1$.
\end{proof}

By \eqref{eq:pyrbounds}, \eqref{eq:general} and \eqref{eq:Lbig},
\begin{equation}\label{eq:chain}
\pyr(d)\ \le\ g(d-1)+1\ \le\ \binom{d+1}{2}+1\ \le\ M(d)+1\ \le\ L(d).
\end{equation}
Hence Theorem~\ref{thm:formula} gives $L(d)\le\iso(d)\le L(d)+1$, and $\iso(d)=L(d)$ unless $\pyr(d)=\binom{d+1}{2}+1$, that is, unless there is a two-distance pyramid in $\R^d$ whose base is a two-distance set with $\binom{d+1}{2}$ points in a hyperplane. Theorem~\ref{thm:nopyr}, which is proved in Section~\ref{sec:nopyr}, shows that there is no such pyramid. For $d\le23$ this can also be seen from the known values of $g$: by \eqref{eq:g}, the equality $g(d-1)=\binom{d+1}{2}$ holds for $d\le23$ only if $d=2$ or $d=9$, and in these two dimensions $\pyr(d)+1\le L(d)$ as well, because $L(2)\ge M(2)+1=6>g(1)+2$ and $\pyr(9)\le g(9)=45$.

\begin{remark}\label{rem:nine}
With the values of $g$ and $M$ in Table~\ref{tab:values}, Theorem~\ref{thm:formula} and the preceding paragraph give the values of $\iso(d)$ for $d\le8$ of Kelly, Croft, Kido and Ionin, and they give the values of $\iso(d)$ in Theorem~\ref{thm:main} and the identity $\iso(23)=\max\{g(23),278\}$ of Theorem~\ref{thm:23} a second time, without the enumeration of Lemma~\ref{lem:B} and without the results of Sections~\ref{sec:tight} and~\ref{sec:nopyr}. In dimension $9$ the value $\iso(9)=46$ needs only $g(9)=45$: by the proof of Theorem~\ref{thm:formula}, an isosceles set with $47$ points in $\R^9$ would consist of a two-distance pyramid with $46$ points and the mirror image of its apex. The results of Section~\ref{sec:nine} are still needed for the uniqueness statement of Theorem~\ref{thm:main} in dimension $9$ and for Theorem~\ref{thm:structure} in dimensions $9$ and $10$.
\end{remark}

\section{Two-distance sets attaining Blokhuis' bound}\label{sec:tight}

By \eqref{eq:chain}, the term $\pyr(d)+1$ in Theorem~\ref{thm:formula} matters only if some two-distance set with $\binom{d+1}{2}$ points in $\R^{d-1}$, the largest number allowed by Blokhuis' bound in \eqref{eq:general}, is the base of a two-distance pyramid. To our knowledge the only two-distance sets known to attain Blokhuis' bound in a dimension $n\ge1$ are three equally spaced points on a line and Lison\v{e}k's set $\LL$ in $\R^8$ \cite{Lis97}. In this section we prove some properties of such sets, and in the next section we use them to prove Theorem~\ref{thm:nopyr}.

Throughout this section and the next, $T\subset\R^n$ is a finite set with exactly two squared distances $a<b$. We write $\mathcal P_2$ for the space of real polynomials of degree at most $2$ on $\R^n$, so that $\dim\mathcal P_2=\binom{n+2}{2}$, and for functions $u,v\colon T\to\R$ we put $\langle u,v\rangle_T=\sum_{x\in T}u(x)v(x)$. Let $\mathcal A^{\perp}$ be the space of all $u\colon T\to\R$ that are orthogonal to the restrictions of all polynomials of degree at most $1$, that is,
\[
\sum_{x\in T}u(x)=0\qquad\text{and}\qquad\sum_{x\in T}u(x)\,x=0 ,
\]
and let $h\colon T\to\R$ be the function $h(x)=\nr{x}$. The starting point is the following identity.

\begin{lemma}\label{lem:F1}
For all $u,v\in\mathcal A^{\perp}$,
\[
\sum_{x,y\in T}u(x)v(y)\ip{x}{y}^2=\frac{ab}{4}\,\langle u,v\rangle_T-\frac12\,\langle u,h\rangle_T\,\langle v,h\rangle_T .
\]
\end{lemma}
\begin{proof}
Put $D(x,y)=\nr{x-y}$, so that $2\ip{x}{y}=h(x)+h(y)-D(x,y)$, and
\[
4\ip{x}{y}^2=(h(x)+h(y))^2-2(h(x)+h(y))D(x,y)+D(x,y)^2 .
\]
We multiply by $u(x)v(y)$ and sum over $x,y\in T$. Since $\sum_xu(x)=\sum_yv(y)=0$ and $\sum_xu(x)x=\sum_yv(y)y=0$, we have
\begin{align*}
\sum_{x,y}u(x)v(y)\,(h(x)+h(y))^2&=2\langle u,h\rangle_T\langle v,h\rangle_T,\\
\sum_{x,y}u(x)v(y)\,(h(x)+h(y))\ip{x}{y}&=0,\\
\sum_{x,y}u(x)v(y)\,D(x,y)&=0 .
\end{align*}
So the first term contributes $2\langle u,h\rangle_T\langle v,h\rangle_T$, and the second term, in which we write $D(x,y)=h(x)+h(y)-2\ip{x}{y}$, contributes $-4\langle u,h\rangle_T\langle v,h\rangle_T$. For the third term note that $(D(x,y)-a)(D(x,y)-b)$ equals $ab$ if $x=y$ and $0$ otherwise, so that
\[
D(x,y)^2=(a+b)D(x,y)-ab+ab\,\delta_{xy}\qquad(x,y\in T),
\]
and the third term contributes $ab\langle u,v\rangle_T$. Altogether the sum of the numbers $4u(x)v(y)\ip{x}{y}^2$ is $ab\langle u,v\rangle_T-2\langle u,h\rangle_T\langle v,h\rangle_T$.
\end{proof}

\begin{lemma}\label{lem:interp}
Every function $T\to\R$ is the restriction to $T$ of a polynomial in $\mathcal P_2$. Hence $|T|\le\binom{n+2}{2}$, and if $|T|=\binom{n+2}{2}$, then every function $T\to\R$ is the restriction of exactly one polynomial in $\mathcal P_2$; in particular, a polynomial in $\mathcal P_2$ that vanishes on $T$ is zero.
\end{lemma}
\begin{proof}
Let $v\colon T\to\R$ be orthogonal to the restrictions of all polynomials in $\mathcal P_2$. Then $v\in\mathcal A^{\perp}$, $\langle v,h\rangle_T=0$, and $\sum_{y}v(y)\ip{x}{y}^2=0$ for every $x$. So Lemma~\ref{lem:F1} with $u=v$ gives $0=\frac{ab}{4}\langle v,v\rangle_T$, and $v=0$. Thus the restriction map $\mathcal P_2\to\R^T$ is surjective, and the remaining statements follow by comparing dimensions.
\end{proof}

The inequality $|T|\le\binom{n+2}{2}$ is Blokhuis' bound \cite{Blo84}, so Lemma~\ref{lem:interp} gives another short proof of it, and it shows what is special about the sets that attain the bound: they are unisolvent for interpolation by polynomials of degree at most $2$.

\begin{lemma}\label{lem:tightg}
Suppose that $|T|=\binom{n+2}{2}$, and let $g$ be the orthogonal projection of $h$ to $\mathcal A^{\perp}$. Then
\[
\sum_{y\in T}g(y)\,yy^{\top}=\frac{ab}{2(n+2)}\,I\qquad\text{and}\qquad\langle g,g\rangle_T=\frac{n\,ab}{2(n+2)} .
\]
\end{lemma}
\begin{proof}
Let $G=\sum_{y}g(y)yy^{\top}$ and $\mu=\frac{ab}{4}-\frac12\langle g,g\rangle_T$. For $u\in\mathcal A^{\perp}$ we have $\langle u,h\rangle_T=\langle u,g\rangle_T$, and $\langle g,h\rangle_T=\langle g,g\rangle_T$. So Lemma~\ref{lem:F1} with $v=g$ gives
\[
\sum_{x\in T}u(x)\,x^{\top}Gx=\frac{ab}{4}\langle u,g\rangle_T-\frac12\langle u,g\rangle_T\langle g,g\rangle_T=\mu\sum_{x\in T}u(x)\nr{x}
\]
for all $u\in\mathcal A^{\perp}$.
Hence the restriction to $T$ of the quadratic form $x^{\top}(G-\mu I)x$ is orthogonal to $\mathcal A^{\perp}$, so it coincides on $T$ with a polynomial $\ell$ of degree at most $1$. The polynomial $x^{\top}(G-\mu I)x-\ell(x)$ vanishes on $T$, so it is zero by Lemma~\ref{lem:interp}, and therefore $G=\mu I$. Taking traces, $\langle g,g\rangle_T=\langle g,h\rangle_T=\tr G=n\mu=n\bigl(\frac{ab}{4}-\frac12\langle g,g\rangle_T\bigr)$, which gives the value of $\langle g,g\rangle_T$, and then $\mu=\frac{ab}{2(n+2)}$.
\end{proof}

The function $g$ does not depend on the choice of the origin, because a translation changes $h$ by the restriction of a polynomial of degree at most $1$. The number $\langle g,g\rangle_T$ is the squared distance, in $\R^T$, from $h$ to the restrictions of polynomials of degree at most $1$; it vanishes if and only if $T$ lies on a sphere. For $\LL$, with squared distances $2$ and $4$, it is $\frac{16}{5}$.

Next we restrict the ratio of the two distances. The following inequality is the relative bound of \cite[Proposition~4.1]{BLSTYY}, in a form that does not require $\gamma^2>n+1$. We include the short proof, because we need an intermediate step.

\begin{lemma}\label{lem:relative}
Let $X\subset\R^n$ be a set with $N$ points and exactly two squared distances $a<b$, and put $\gamma=(a+b)/(b-a)$. Then
\[
\gamma^2\,(N-1-n)\ \le\ (n+1)(N-1),
\]
and if equality holds, then $X$ lies on a sphere.
\end{lemma}
\begin{proof}
Scale $X$ so that $b-a=2$; then $a=\gamma-1$ and $b=\gamma+1$. Let $D$ be the matrix of squared distances of $X$, let $\mathbf 1$ be the all-ones vector of length $N$, and let $A$ be the symmetric matrix with zero diagonal whose off-diagonal entries are $+1$ for the pairs at squared distance $b$ and $-1$ for the pairs at squared distance $a$. Then $D=\gamma(\mathbf 1\mathbf 1^{\top}-I)+A$. With $P=I-\frac1N\mathbf 1\mathbf 1^{\top}$, the matrix $G=-\frac12PDP=\frac12(\gamma P-PAP)$ is the Gram matrix of $X$ with respect to its centroid, so it is positive semidefinite of rank at most $n$, and therefore $(\tr G)^2\le n\,\tr(G^2)$. Put $\sigma=\frac1N\mathbf 1^{\top}A\mathbf 1$. Using $\tr A=0$ and $\tr A^2=N(N-1)$ we get $\tr(PAP)=-\sigma$ and $\tr\bigl((PAP)^2\bigr)=N(N-1)-\frac2N\lVert A\mathbf 1\rVert^2+\sigma^2$, hence $2\tr G=\gamma(N-1)+\sigma$ and
\begin{equation}\label{eq:trace}
(\gamma(N-1)+\sigma)^2\ \le\ n\Bigl(\gamma^2(N-1)+2\gamma\sigma+N(N-1)-\tfrac2N\lVert A\mathbf 1\rVert^2+\sigma^2\Bigr).
\end{equation}
By the Cauchy--Schwarz inequality $\lVert A\mathbf 1\rVert^2\ge N\sigma^2$, so
\[
(n+1)\sigma^2+2\gamma(N-1-n)\sigma+\gamma^2(N-1)(N-1-n)-nN(N-1)\ \le\ 0 .
\]
The minimum over $\sigma$ of the left-hand side is $\gamma^2(N-1-n)\frac{nN}{n+1}-nN(N-1)$, and the inequality follows. If equality holds, then $\lVert A\mathbf 1\rVert^2=N\sigma^2$, so $A\mathbf 1$ and hence $D\mathbf 1$ is a multiple of $\mathbf 1$. Since
\begin{equation}\label{eq:centroid}
\sum_{y\in X}\nr{x-y}=N\nr{x-c}+\sum_{y\in X}\nr{y-c}
\end{equation}
for the centroid $c$ of $X$, all points of $X$ are then at the same distance from $c$.
\end{proof}

For a set with $N=\binom{n+2}{2}$ points, Lemma~\ref{lem:relative} gives $\gamma^2\le n+3$, because then $N-1-n=\frac{n(n+1)}{2}$ and $N-1=\frac{n(n+3)}{2}$. Lemma~\ref{lem:tightg} allows us to improve this.

\begin{proposition}\label{prop:tight}
Let $T\subset\R^n$ be a two-distance set with $\binom{n+2}{2}$ points and squared distances $a<b$, and put $\gamma=(a+b)/(b-a)$. Then
\begin{equation}\label{eq:tightgamma}
\gamma^2\,(n+4)\ \le\ n^2+5n+8 .
\end{equation}
If $n\ge9$, then $a/b=(\kappa-1)/\kappa$ for an integer $\kappa\ge3$ with $(2\kappa-1)^2\le n+1$; in particular $n\ge24$.
\end{proposition}
\begin{proof}
Put $N=\binom{n+2}{2}$, scale $T$ so that $b-a=2$, and use the notation of the proof of Lemma~\ref{lem:relative} with $X=T$. Take the centroid of $T$ as the origin, so that $h(x)$ is the squared distance of $x$ from the centroid. By \eqref{eq:centroid} and $D=\gamma(\mathbf 1\mathbf 1^{\top}-I)+A$, the function $h-\frac1NA\mathbf 1$ is constant on $T$. Constant functions are orthogonal to $\mathcal A^{\perp}$, so the function $g$ of Lemma~\ref{lem:tightg} is also the orthogonal projection of $\frac1N(A\mathbf 1-\sigma\mathbf 1)$ to $\mathcal A^{\perp}$, and therefore
\[
N^2\langle g,g\rangle_T\ \le\ \lVert A\mathbf 1-\sigma\mathbf 1\rVert^2=\lVert A\mathbf 1\rVert^2-N\sigma^2 .
\]
Since $ab=\gamma^2-1$, Lemma~\ref{lem:tightg} gives $\lVert A\mathbf 1\rVert^2\ge N\sigma^2+N^2\frac{n(\gamma^2-1)}{2(n+2)}$. Inserting this in \eqref{eq:trace} we get
\begin{multline*}
(n+1)\sigma^2+2\gamma(N-1-n)\sigma+\gamma^2(N-1)(N-1-n)\\
-nN(N-1)+\frac{n^2N(\gamma^2-1)}{n+2}\ \le\ 0 .
\end{multline*}
The sum of the first two terms is at least $-\gamma^2(N-1-n)^2/(n+1)$. With $N-1-n=\frac{n(n+1)}{2}$, $N-1=\frac{n(n+3)}{2}$ and $N=\frac{(n+1)(n+2)}{2}$, division by $\frac{n^2(n+1)}{4}$ gives
\[
-\gamma^2+\gamma^2(n+3)-(n+2)(n+3)+2(\gamma^2-1)\ \le\ 0,
\]
which is \eqref{eq:tightgamma}.

Now let $n\ge9$. Since $N>2n+3$, the theorem of Larman, Rogers and Seidel \cite{LRS77} shows that $a/b=(\kappa-1)/\kappa$ for an integer $\kappa\ge2$; thus $\gamma=2\kappa-1$. If $\kappa=2$, the squared distances are in ratio $1:2$, and then $|T|\le\binom{n+1}{2}<N$ by \cite[Theorem~1.1]{BLSTYY}, as $n\ge9$. Hence $\kappa\ge3$. By \eqref{eq:tightgamma}, $\gamma^2\le n+1+\frac{4}{n+4}<n+2$, and as $\gamma$ is an integer, $(2\kappa-1)^2\le n+1$. With $\kappa\ge3$ this gives $n\ge24$.
\end{proof}

For three equally spaced points on a line, \eqref{eq:tightgamma} reads $\frac{125}{9}\le14$, and for $\LL$, where $\gamma=3$ and $n=8$, it reads $108\le112$. So \eqref{eq:tightgamma} is close to an equality for both known examples.

\begin{corollary}\label{cor:g23}
For $9\le n\le23$ there is no two-distance set with $\binom{n+2}{2}$ points in $\R^n$. In particular $277\le g(23)\le299$.
\end{corollary}
\begin{proof}
The first statement is contained in Proposition~\ref{prop:tight}. It gives $g(23)\le\binom{25}{2}-1=299$ by Blokhuis' bound, and $g(23)\ge277$ by \cite{GKM26}.
\end{proof}

For $n\le22$ Corollary~\ref{cor:g23} is contained in \eqref{eq:g}. By Proposition~\ref{prop:tight}, a two-distance set attaining Blokhuis' bound in a dimension $n\ge9$ has squared distances in ratio $2:3$ and $n\ge24$, or in ratio $3:4$ and $n\ge48$, and so on. We do not know whether such sets exist; the first open case is a two-distance set with $325$ points in $\R^{24}$.

Finally we consider sets attaining Blokhuis' bound that lie on two concentric spheres, as $\LL$ does. Let $|T|=\binom{n+2}{2}$, and suppose that $T$ is contained in the union of two spheres with a common centre $m$ and squared radii $r_1\ne r_2$, where a sphere of radius zero is the point $m$. Put
\[
T_i=\{x\in T:\nr{x-m}=r_i\}\qquad(i=1,2).
\]
A two-distance set on a sphere in $\R^n$ has at most $\frac{n(n+3)}{2}<\binom{n+2}{2}$ points by \eqref{eq:general}, so $T_1$ and $T_2$ are both nonempty. For real $r$ let
\[
\varphi(r)=(2r-a)(2r-b)+\frac{4r^2}{n} .
\]

\begin{theorem}\label{thm:twosphere}
In this situation suppose that $r_1>0$, $r_2>0$ and $r_1+r_2\ne\frac{a+b}{2}$. Then:
\begin{enumerate}[label=\textup{(\roman*)}]
\item $4r_1r_2=n\,(r_1+r_2-a)(b-r_1-r_2)$;
\item for $i=1,2$ the set $T_i$ is a spherical $2$-design, that is,
\[
\sum_{y\in T_i}(y-m)=0\qquad\text{and}\qquad\sum_{y\in T_i}(y-m)(y-m)^{\top}=\frac{r_i\,|T_i|}{n}\,I ;
\]
in particular $m$ is the centroid of $T_1$, of $T_2$ and of $T$;
\item $|T_i|\,\varphi(r_i)=ab$ for $i=1,2$.
\end{enumerate}
\end{theorem}
\begin{proof}
We may assume that $m=0$. Put $\lambda=r_1+r_2$. By symmetry it suffices to prove (i), and (ii) and (iii) for $i=1$.

Let $Z\subseteq\mathcal P_2$ be the space of the polynomials in $\mathcal P_2$ that vanish on $T_1$. By Lemma~\ref{lem:interp} the restriction map $\mathcal P_2\to\R^T$ is bijective. Hence the restriction map $\mathcal P_2\to\R^{T_1}$ is surjective with kernel $Z$, and $\dim Z=|T_2|$. For $y\in T_2$ let
\begin{align*}
h_y(x)&=\bigl(\lambda-a-2\ip{x}{y}\bigr)\bigl(\lambda-b-2\ip{x}{y}\bigr)\\
&=4\ip{x}{y}^2-2(2\lambda-a-b)\ip{x}{y}+(\lambda-a)(\lambda-b).
\end{align*}
If $x\in T_1$, then $\nr{x-y}=\lambda-2\ip{x}{y}$ is $a$ or $b$, so $h_y(x)=0$. Thus $h_y\in Z$.

\emph{Step 1: the polynomials $h_y$, $y\in T_2$, form a basis of $Z$.} As $\dim Z=|T_2|$, we have to show that they are linearly independent. Let $\sum_{y\in T_2}c_yh_y=0$ with real numbers $c_y$. Comparing the terms of degree $2$ and of degree $1$, and using $2\lambda\ne a+b$, we get
\begin{equation}\label{eq:dep}
\sum_{y\in T_2}c_y\ip{x}{y}^2=0\quad(x\in\R^n)\qquad\text{and}\qquad\sum_{y\in T_2}c_y\,y=0 .
\end{equation}
For $z\in T_2$ let $F_z(x)=(2r_2-a-2\ip{x}{z})(2r_2-b-2\ip{x}{z})$. If $x\in T_2$, then $\nr{x-z}=2r_2-2\ip{x}{z}$, so $F_z(x)=ab$ if $x=z$ and $F_z(x)=0$ otherwise. On the other hand, by \eqref{eq:dep},
\[
\sum_{z\in T_2}c_zF_z(x)=(2r_2-a)(2r_2-b)\sum_{z\in T_2}c_z\qquad(x\in\R^n).
\]
Evaluating at the points $x\in T_2$ we see that $ab\,c_x$ does not depend on $x$. So all $c_y$ are equal to one number $c$, and \eqref{eq:dep} gives $c\sum_{y\in T_2}\ip{x}{y}^2=0$ for all $x$. Taking $x\in T_2$, where $\ip{x}{x}^2=r_2^2>0$, we get $c=0$.

\emph{Step 2: proof of \textup{(i)}.} The polynomial $\nr{x}-r_1$ belongs to $Z$, so by Step~1 there are real numbers $c_y$ with $\nr{x}-r_1=\sum_{y\in T_2}c_yh_y(x)$. Applying the Laplace operator $\Delta$ to both sides gives $2n=8\sum_yc_y\nr{y}=8r_2\sum_yc_y$, and comparing the constant terms gives $-r_1=(\lambda-a)(\lambda-b)\sum_yc_y$. Hence $-4r_1r_2=n(\lambda-a)(\lambda-b)$, which is (i).

\emph{Step 3: proof of \textup{(ii)} and \textup{(iii)}.} Consider the linear functional
\[
\Lambda(f)=f(0)+\frac{r_1}{2n}\,\Delta f(0)\qquad(f\in\mathcal P_2).
\]
(It is the mean value of $f$ over the sphere $\nr{x}=r_1$, but we do not need this.) For $y\in T_2$ we have $\Lambda(h_y)=(\lambda-a)(\lambda-b)+\frac{r_1}{2n}\cdot8r_2=0$ by (i). So $\Lambda$ vanishes on $Z$ by Step~1, and since $Z$ is the kernel of the surjective restriction map $\mathcal P_2\to\R^{T_1}$, there is a function $\omega\colon T_1\to\R$ with
\begin{equation}\label{eq:cubature}
\Lambda(f)=\sum_{y\in T_1}\omega(y)f(y)\qquad(f\in\mathcal P_2).
\end{equation}
For $z\in T_1$ let $G_z(x)=(2r_1-a-2\ip{x}{z})(2r_1-b-2\ip{x}{z})$. As in Step~1, $G_z(x)=ab$ if $x=z$ and $G_z(x)=0$ for the other points $x\in T_1$, so \eqref{eq:cubature} gives $\Lambda(G_z)=ab\,\omega(z)$. On the other hand $G_z(0)=(2r_1-a)(2r_1-b)$ and $\Delta G_z=8\nr{z}=8r_1$, so $\Lambda(G_z)=\varphi(r_1)$. Hence $\omega(z)=\varphi(r_1)/ab$ for every $z\in T_1$, and \eqref{eq:cubature} becomes
\[
ab\,\Lambda(f)=\varphi(r_1)\sum_{y\in T_1}f(y)\qquad(f\in\mathcal P_2).
\]
For $f=1$ this is (iii). For $f(x)=x_j$ and for $f(x)=x_jx_k$ we have $\Lambda(f)=0$ and $\Lambda(f)=\frac{r_1}{n}\delta_{jk}$, respectively, and (ii) follows by (iii).
\end{proof}

\begin{example}\label{ex:lisonek}
Up to similarity, $\LL$ consists of the $36$ points $e_i+e_j$ $(1\le i<j\le9)$ and the $9$ points $\frac13\mathbf 1-e_i$ $(1\le i\le9)$ in the hyperplane $\sum_ix_i=2$ of $\R^9$, where $e_1,\dots,e_9$ is the standard basis and $\mathbf 1=\sum_ie_i$; here $n=8$, $a=2$ and $b=4$. Its points lie on two spheres about the centroid $\frac29\mathbf 1$, with $r_1=\frac{14}{9}$ for the first $36$ points and $r_2=\frac89$ for the other $9$. Both sides of (i) are $\frac{448}{81}$, and $\varphi(\frac{14}{9})=\frac29=\frac{8}{36}$ and $\varphi(\frac89)=\frac89$, in accordance with (iii). The assumption $r_1+r_2\ne\frac{a+b}{2}$ cannot be omitted: for $T=\{0,1,2\}\subset\R$, with $a=1$ and $b=4$, and $m=\frac12$ we have $T_1=\{0,1\}$, $T_2=\{2\}$, $r_1=\frac14$, $r_2=\frac94$ and $r_1+r_2=\frac{a+b}{2}$, but $m$ is not the centroid of $T_2$, and $|T_2|\varphi(r_2)=22\ne ab$.
\end{example}

Thus the radii of the two spheres determine the numbers of points on them, and by (i) and (iii) the two radii satisfy the two equations
\[
4r_1r_2=n(r_1+r_2-a)(b-r_1-r_2)\qquad\text{and}\qquad\frac{ab}{\varphi(r_1)}+\frac{ab}{\varphi(r_2)}=\binom{n+2}{2}.
\]
In the next section we show that these equations have no solution with $r_2-r_1=b-a$.

\section{Proof of Theorems~\ref{thm:all} and~\ref{thm:nopyr}}\label{sec:nopyr}

We keep the notation of Section~\ref{sec:tight}. Besides Theorem~\ref{thm:twosphere} and Lemma~\ref{lem:relative} we need the following bound of Musin \cite{Mus09} for spherical two-distance sets whose inner products have a nonnegative sum. We recall its short proof.

\begin{lemma}[Musin]\label{lem:musin}
Let $r>0$, let $S$ be a finite subset of the sphere $\{x\in\R^n:\nr{x}=r\}$, and let $\alpha,\beta<r$ be real numbers with $\alpha+\beta\ge0$ such that $\ip{x}{y}\in\{\alpha,\beta\}$ for all distinct $x,y\in S$. Then $|S|\le\binom{n+1}{2}$.
\end{lemma}
\begin{proof}
Let $\Sigma$ be the sphere $\nr{x}=r$ and let $W$ be the space of all functions on $\Sigma$ of the form $x\mapsto q(x)+\ip{u}{x}$ with a quadratic form $q$ and a vector $u\in\R^n$; then $\dim W\le\binom{n+1}{2}+n$. For $z\in S$ let $f_z(x)=(\ip{x}{z}-\alpha)(\ip{x}{z}-\beta)$. The restriction of $f_z$ to $\Sigma$ belongs to $W$, because on $\Sigma$ the constant $\alpha\beta$ coincides with the quadratic form $\frac{\alpha\beta}{r}\nr{x}$. For $x\in S$ we have $f_z(x)=\theta\,\delta_{xz}$ with $\theta=(r-\alpha)(r-\beta)>0$. We claim that the functions $f_z$ $(z\in S)$ and the $n$ coordinate functions are linearly independent on $\Sigma$; then $|S|+n\le\binom{n+1}{2}+n$. Suppose that $\sum_{z\in S}c_zf_z(x)+\ip{u}{x}=0$ for all $x\in\Sigma$. Replacing $x$ by $-x$ and subtracting, we get $-2(\alpha+\beta)\sum_zc_z\ip{x}{z}+2\ip{u}{x}=0$ for all $x\in\Sigma$, so $u=(\alpha+\beta)\sum_zc_zz$. Evaluating at $x\in S$ we get $\theta c_x=-\ip{u}{x}$. Hence
\[
\nr{u}=(\alpha+\beta)\sum_{z\in S}c_z\ip{u}{z}=-\frac{\alpha+\beta}{\theta}\sum_{z\in S}\ip{u}{z}^2\ \le\ 0,
\]
so $u=0$, and then $c_x=0$ for all $x\in S$.
\end{proof}

\begin{proof}[Proof of Theorem~\ref{thm:nopyr}]
Put $N=\binom{n+2}{2}$. Suppose that there are $m\in\R^n$ and $t\in\R$ with $\nr{x-m}+t\in\{a,b\}$ for all $x\in T$; we may assume that $m=0$. Put $r_1=a-t$ and $r_2=b-t$, so that $r_2-r_1=b-a>0$, and $T_i=\{x\in T:\nr{x}=r_i\}$. Then $T=T_1\cup T_2$. If $T_1$ or $T_2$ were empty, then $T$ would lie on a sphere, which is impossible because a two-distance set on a sphere in $\R^n$ has at most $\frac{n(n+3)}{2}<N$ points by \eqref{eq:general}. So $T_1$ and $T_2$ are nonempty; in particular $r_1\ge0$. For distinct points $x,y\in T_i$ we have
\begin{equation}\label{eq:inner}
\ip{x}{y}=r_i-\tfrac12\nr{x-y}\in\bigl\{r_i-\tfrac a2,\ r_i-\tfrac b2\bigr\}.
\end{equation}
We distinguish three cases.

\emph{Case 1: $r_1=0$.} Then $T_1=\{0\}$, and every $x\in T_2$ satisfies $\nr{x-0}=r_2=b-a$. As $0\in T$, this number is $a$ or $b$, and it is smaller than $b$, so $b=2a$ and $r_2=a$. By \eqref{eq:inner} the inner products of distinct points of $T_2$ are $\frac a2$ and $0$, so $|T_2|\le\binom{n+1}{2}$ by Lemma~\ref{lem:musin}. But $|T_2|=N-1=\binom{n+1}{2}+n$, a contradiction.

\emph{Case 2: $r_1>0$ and $r_1+r_2=\frac{a+b}{2}$.} Then $r_1=\frac{3a-b}{4}$ and $r_2=\frac{3b-a}{4}$. By \eqref{eq:inner} the inner products of distinct points of $T_1$ are $-\frac{b-a}{4}$ and $-\frac{3(b-a)}{4}$, and those of distinct points of $T_2$ are $\frac{3(b-a)}{4}$ and $\frac{b-a}{4}$, which are positive and smaller than $r_2$. So $|T_2|\le\binom{n+1}{2}$ by Lemma~\ref{lem:musin}, and
\[
0\ \le\ \Bigl\lVert\sum_{x\in T_1}x\Bigr\rVert^2\ \le\ |T_1|\,r_1-|T_1|(|T_1|-1)\frac{b-a}{4},
\]
so $|T_1|\le1+\frac{4r_1}{b-a}=\gamma-1$, where $\gamma=\frac{a+b}{b-a}$. Hence $N=|T_1|+|T_2|\le\gamma-1+\binom{n+1}{2}$, that is, $\gamma\ge n+2$. But $\gamma^2\le n+3$ by Lemma~\ref{lem:relative}, as noted after its proof, and $(n+2)^2>n+3$, a contradiction.

\emph{Case 3: $r_1>0$ and $r_1+r_2\ne\frac{a+b}{2}$.} Then Theorem~\ref{thm:twosphere} applies. As $r_1+r_2-a=b-2t$ and $b-r_1-r_2=2t-a$, part (i) reads $4(a-t)(b-t)=n(b-2t)(2t-a)$, that is,
\begin{equation}\label{eq:tquad}
4(n+1)t^2-2(n+2)(a+b)t+(n+4)ab=0 .
\end{equation}
Put $U=4(n+1)t-2(n+2)a$ and $V=2(n+2)b-4(n+1)t$, so that $U+V=2(n+2)(b-a)$. We claim that
\begin{equation}\label{eq:UV}
n\,\varphi(a-t)=(b-a)\,U,\qquad n\,\varphi(b-t)=(b-a)\,V,\qquad UV=4n\,ab .
\end{equation}
Indeed,
\begin{align*}
n\,\varphi(a-t)&=n(a-2t)(2a-b-2t)+4(a-t)^2\\
&=4(n+1)t^2-(6na-2nb+8a)t+(2n+4)a^2-nab,
\end{align*}
and replacing $4(n+1)t^2$ by $2(n+2)(a+b)t-(n+4)ab$ according to \eqref{eq:tquad} gives $4(n+1)(b-a)t-2(n+2)a(b-a)=(b-a)U$. The second identity follows by exchanging $a$ and $b$, which leaves $\varphi$ and \eqref{eq:tquad} unchanged. Finally, by \eqref{eq:tquad},
\begin{align*}
UV&=-4(n+1)\bigl(4(n+1)t^2-2(n+2)(a+b)t\bigr)-4(n+2)^2ab\\
&=4(n+1)(n+4)ab-4(n+2)^2ab=4n\,ab .
\end{align*}
By Theorem~\ref{thm:twosphere}(iii) and \eqref{eq:UV},
\[
N=|T_1|+|T_2|=\frac{ab}{\varphi(a-t)}+\frac{ab}{\varphi(b-t)}=\frac{n\,ab}{b-a}\cdot\frac{U+V}{UV}=\frac{n+2}{2}.
\]
This contradicts $N=\frac{(n+1)(n+2)}{2}$, because $n\ge1$. This proves the first statement of the theorem.

Now let $X\subset\R^{n+1}$ be a two-distance pyramid with apex $p$ whose base $T$ lies in a hyperplane $H$ and has $\binom{n+2}{2}$ points. We identify $H$ with $\R^n$. As $T$ has more than $n+1$ points, it has exactly two squared distances $a<b$, and these are the squared distances of $X$. Let $m$ be the orthogonal projection of $p$ to $H$ and $t=\nr{p-m}$. Then $\nr{x-m}+t=\nr{x-p}\in\{a,b\}$ for all $x\in T$, which is impossible by the first statement. Finally let $d\ge2$. By \eqref{eq:pyrbounds}, $\pyr(d)\le\binom{d+1}{2}+1$, and a two-distance pyramid with $\binom{d+1}{2}+1$ points in $\R^d$ would have a base with $\binom{d+1}{2}=\binom{(d-1)+2}{2}$ points in a hyperplane. Hence $\pyr(d)\le\binom{d+1}{2}$.
\end{proof}

\begin{proof}[Proof of Theorem~\ref{thm:all}]
Let $L(d)$ be the right-hand side of \eqref{eq:forum}. By Theorem~\ref{thm:formula}, $\iso(d)=\max\{L(d),\pyr(d)+1\}$, and $\pyr(d)+1\le\binom{d+1}{2}+1\le L(d)$ by Theorem~\ref{thm:nopyr} and \eqref{eq:Lbig}.
\end{proof}

The proof gives slightly more: the isosceles sets that are not two-distance sets satisfy the bound without the term $g(d)$.

\begin{corollary}\label{cor:three}
Let $d\ge2$. Every isosceles set in $\R^d$ that is not a two-distance set has at most $\max\{M(d)+1,\ M(d-1)+3\}$ points.
\end{corollary}
\begin{proof}
Put $L'(d)=\max\{M(d)+1,M(d-1)+3\}$; then $L'(d)\ge\binom{d+1}{2}+1$ by \eqref{eq:general}. Let $S\subset\R^d$ be an isosceles set that is not a two-distance set, and let $(S_1,\dots,S_k)$ be a complete decomposition of $S$, of type $a$. Then $k\ge2$, so $a\ne(d)$. If $a$ is not special, then $|S|\le v(a)\le\binom{d+1}{2}+1\le L'(d)$ by Lemma~\ref{lem:nonspecial}. If $a$ is special and $a\ne(1,d-1)$, then $|S|\le L'(d)$ by the proof of Lemma~\ref{lem:types}. Finally let $a=(1,d-1)$, and argue as in the proof of Theorem~\ref{thm:formula}. In case (a) of Lemma~\ref{lem:twopoints}, $|S|\le M(d-1)+2$. In cases (b) and (c), $S$ has a complete decomposition with at least two parts whose type is not $(1,d-1)$, so $|S|\le L'(d)$ by what we have already shown. In case (d), $|S|\le\pyr(d)+1\le\binom{d+1}{2}+1\le L'(d)$ by Theorem~\ref{thm:nopyr}.
\end{proof}

\begin{remark}\label{rem:cases}
Case~3 is the general case, and there the contradiction is quantitative: the base of a two-distance pyramid would have to consist of two spherical $2$-designs with $\frac{n+2}{2}$ points altogether. Cases~1 and~2 correspond to the two exceptions in Theorem~\ref{thm:twosphere}.
\end{remark}

\begin{proof}[Proof of Corollary~\ref{cor:GY}]
By the theorem of Glazyrin and Yu \cite{GY18} quoted in Section~\ref{sec:inputs}, $M(d)=\binom{d+1}{2}$ and $M(d-1)=\binom d2$. So $M(d-1)+3\le M(d)+1$, and the claim follows from Theorem~\ref{thm:all}.
\end{proof}

\section{Large isosceles sets in every dimension}\label{sec:large}

Theorem~\ref{thm:structure} describes the large isosceles sets in dimensions $10\le d\le22$, where $g$ and $M$ are known. The arguments of Sections~\ref{sec:formula} and~\ref{sec:nopyr} give a weaker description that uses only the general bounds \eqref{eq:general} and is valid in every dimension.

\begin{theorem}\label{thm:large}
Let $d\ge5$ and let $S\subset\R^d$ be an isosceles set with $|S|\ge\binom{d+1}{2}+1$. Then at least one of the following holds.
\begin{enumerate}[label=\textup{(\roman*)}]
\item $S$ is a two-distance set.
\item $S=X\cup\{o\}$, where $X$ is a two-distance set of dimension $d$ on a sphere with centre $o$.
\item $S=X\cup Y$, where $X$ is a two-distance set of dimension $d-1$ lying on a sphere and $|Y|\le3$. This requires $M(d-1)\ge\binom{d+1}{2}-2$.
\item $|S|=\binom{d+1}{2}+1$ and $S=T\cup\{p,q\}$, where $T$ is a two-distance set with $\binom{d+1}{2}-1$ points that spans a hyperplane $H$, $T\cup\{p\}$ is a two-distance set, $p\notin H$, and $q$ is the mirror image of $p$ in $H$. This requires $g(d-1)\ge\binom{d+1}{2}-1$.
\end{enumerate}
In cases \textup{(i)}, \textup{(ii)} and \textup{(iii)}, $|S|$ is at most $g(d)$, $M(d)+1$ and $M(d-1)+3$, respectively.
\end{theorem}
\begin{proof}
Since $|S|>\binom{d+1}{2}$ and $d\ge5$, every complete decomposition of $S$ has a special type, by Lemma~\ref{lem:types} for the types different from $(1,d-1)$. Let $(S_1,\dots,S_k)$ be a complete decomposition of $S$, of type $a$. Recall that for $i<k$ the set $S_i$ lies on a sphere centred at a point of $S_k$, that $|S_i|\le2$ if $\dim S_i=1$ and $i<k$, that $|S_k|\le3$ if $\dim S_k=1$, and that $|S_k|=1$ if $\dim S_k=0$.

If $a=(d)$, then (i) holds. If $a=(d,0)$, then (ii) holds with $X=S_1$ and $S_2=\{o\}$. If $a=(d-1,1)$, then (iii) holds with $X=S_1$ and $Y=S_2$. If $a=(d-1,1,0)$, then (iii) holds with $X=S_1$ and $Y=S_2\cup S_3$. If $a=(1,d-1,0)$, then (iii) holds with $X=S_2$ and $Y=S_1\cup S_3$. In case (iii), $|X|\ge|S|-3\ge\binom{d+1}{2}-2$.

Let $a=(1,d-1)$. As in the proof of Theorem~\ref{thm:formula}, $S=T\cup\{p,q\}$ with $T=S_2$ spanning the bisecting hyperplane $H$ of $p$ and $q$, and we apply Lemma~\ref{lem:twopoints}. In case (a), (iii) holds with $X=T$ and $Y=\{p,q\}$. In case (b) the complete decomposition $(\{p,q\},C,T\setminus C)$ has a special type with three parts and first part of dimension $1$, which must be $(1,d-1,0)$, and (iii) holds as above. Case (c) cannot occur, because a type whose first entry is $2$ is not special when $d\ge5$. In case (d), $T\cup\{p\}$ is a two-distance pyramid, so $|T|+1\le\pyr(d)\le\binom{d+1}{2}$ by Theorem~\ref{thm:nopyr}. As $|S|=|T|+2\ge\binom{d+1}{2}+1$, this gives $|T|=\binom{d+1}{2}-1$ and $|S|=\binom{d+1}{2}+1$, and (iv) holds. The bounds for $|S|$ in cases (i)--(iii) follow from \eqref{eq:parts}.
\end{proof}

\begin{corollary}\label{cor:unique}
Let $d\ge5$ and suppose that $g(d)=M(d)=\binom{d+1}{2}$, $M(d-1)\le\binom{d+1}{2}-3$ and $g(d-1)\le\binom{d+1}{2}-2$. Then $\iso(d)=\binom{d+1}{2}+1$, and every isosceles set of this size in $\R^d$ consists of a two-distance set with $\binom{d+1}{2}$ points on a sphere together with the centre of the sphere.
\end{corollary}
\begin{proof}
The hypotheses and Theorem~\ref{thm:all} give $\iso(d)=\binom{d+1}{2}+1$. Let $S\subset\R^d$ be isosceles with $\binom{d+1}{2}+1$ points and apply Theorem~\ref{thm:large}. Case (i) is impossible because $|S|>g(d)$, case (iii) because $M(d-1)+3\le\binom{d+1}{2}$, and case (iv) because $g(d-1)<\binom{d+1}{2}-1$. So (ii) holds, and $|X|=\binom{d+1}{2}$.
\end{proof}

The hypotheses of Corollary~\ref{cor:unique} hold for $10\le d\le21$ by Table~\ref{tab:values}, which gives Theorem~\ref{thm:main}(a) in these dimensions once more, up to the identification of the spherical set with $\J(d+1,2)$. The hypothesis on $g(d-1)$ is needed only to exclude case (iv) of Theorem~\ref{thm:large}; see Problem~\ref{prob:near}.

\section{Open questions}\label{sec:open}

By Theorem~\ref{thm:all}, the problem of Erd\H{o}s to determine $\iso(d)$ for all $d$ is reduced to the determination of $g$ and $M$. What remains open is therefore the following.
\begin{enumerate}[label=\textup{(\arabic*)}]
\item The values of $g(d)$ for $d\ge23$. The first case is $277\le g(23)\le299$ (Corollary~\ref{cor:g23}), and by Theorem~\ref{thm:23} it decides the value of $\iso(23)$.
\item The values of $M(d)$ for $d=(2k+1)^2-3\ge46$, the dimensions left open in \cite{GY18}; they enter Theorem~\ref{thm:all} in the dimensions $d$ and $d+1$.
\end{enumerate}
Theorem~\ref{thm:large} leaves two questions about the largest isosceles sets.

\begin{problem}\label{prob:tight}
Is there, for some $n\ge24$, a two-distance set with $\binom{n+2}{2}$ points in $\R^n$?
\end{problem}

By Proposition~\ref{prop:tight} such a set has squared distances in ratio $(\kappa-1):\kappa$ for an integer $\kappa\ge3$ with $(2\kappa-1)^2\le n+1$, and by Lemma~\ref{lem:interp} it is unisolvent for polynomials of degree at most $2$. If there is no such set, then $g(d)\le\binom{d+2}{2}-1$ for all $d\ge9$.

\begin{problem}\label{prob:near}
Is there, for some $n\ge9$, a two-distance set with $\binom{n+2}{2}-1$ points in $\R^n$ that does not lie on a sphere and is the base of a two-distance pyramid?
\end{problem}

If there is no such set, then in case (iv) of Theorem~\ref{thm:large} the set $T$ lies on a sphere, so that (iii) holds as well, and the hypothesis on $g(d-1)$ in Corollary~\ref{cor:unique} can be omitted. Spherical bases do occur: a two-distance set with $\binom{n+2}{2}-1=\frac{n(n+3)}{2}$ points on a sphere attains the bound of Delsarte, Goethals and Seidel in \eqref{eq:general}, and for example the McLaughlin configuration in $\R^{22}$, whose squared distances are $\frac53$ and $\frac52$ if the sphere has radius $1$, is the base of a two-distance pyramid with $276$ points in $\R^{23}$, with apex at squared distance $\frac53$ from all its points. The proof of Theorem~\ref{thm:nopyr} does not extend directly to sets with $\binom{n+2}{2}-1$ points, because such a set lies on a quadric.


\normalsize

\begin{thebibliography}{22}

\normalsize

\bibitem{BY13} A. Barg and W.-H. Yu,
\emph{New bounds for spherical two-distance sets},
Exp. Math. 22 (2013), 187--194.\newline
DOI: 10.1080/10586458.2013.767725

\bibitem{Blo83} A. Blokhuis,
\emph{Few-distance sets},
Ph.D. thesis, Eindhoven University of Technology, 1983; CWI Tract 7, Centrum voor Wiskunde en Informatica, Amsterdam, 1984.

\bibitem{Blo84} A. Blokhuis,
\emph{A new upper bound for the cardinality of $2$-distance sets in Euclidean space},
in: Convexity and Graph Theory (Jerusalem, 1981), Ann. Discrete Math. 20, North-Holland, Amsterdam, 1984, 65--66.

\bibitem{EP503} T. F. Bloom,
\emph{Erd\H{o}s Problem \#503},
\url{https://www.erdosproblems.com/503}, and its forum thread \url{https://www.erdosproblems.com/forum/thread/503}, accessed November 2026.

\bibitem{BLSTYY} B. Bu, C.-Y. Lee, Y.-H. Sheng, M.-T. Tsai, H. Yi and W.-H. Yu,
\emph{Johnson optimality and uniqueness for distance ratio $1:\sqrt2$},
arXiv:2610.40367v1 (2026).

\bibitem{CY25} W.-C. Chen and W.-H. Yu,
\emph{Bounds on two-distance sets in Euclidean space and unit sphere},
arXiv:2509.00858 (2025).

\bibitem{Cro62} H. T. Croft,
\emph{$9$-point and $7$-point configurations in $3$-space},
Proc. London Math. Soc. (3) 12 (1962), 400--424.

\bibitem{DGS77} P. Delsarte, J. M. Goethals and J. J. Seidel,
\emph{Spherical codes and designs},
Geom. Dedicata 6 (1977), 363--388.\newline
DOI: 10.1007/BF03187604

\bibitem{ES66} S. J. Einhorn and I. J. Schoenberg,
\emph{On Euclidean sets having only two distances between points I, II},
Indag. Math. 28 (1966), 479--488, 489--504.

\bibitem{Erd46} P. Erd\H{o}s,
\emph{Problem E735},
Amer. Math. Monthly 53 (1946), 394.

\bibitem{EK47} P. Erd\H{o}s and L. M. Kelly,
\emph{Isosceles $n$-points (solution of Problem E735)},
Amer. Math. Monthly 54 (1947), 227--229.

\bibitem{Gal67} T. Gallai,
\emph{Transitiv orientierbare Graphen},
Acta Math. Acad. Sci. Hungar. 18 (1967), 25--66; English transl. by F. Maffray and M. Preissmann in: J. L. Ram\'irez Alfons\'in and B. A. Reed (eds.), Perfect Graphs, Wiley, 2001, 25--66.

\bibitem{GKM26} H.-J. Ge, J. Koolen and A. Munemasa,
\emph{A $2$-distance set with $277$ points in the Euclidean space of dimension $23$},
Discrete Comput. Geom. 76 (2026), 1679--1686.\newline
DOI: 10.1007/s00454-026-00843-9

\bibitem{GY18} A. Glazyrin and W.-H. Yu,
\emph{Upper bounds for $s$-distance sets and equiangular lines},
Adv. Math. 330 (2018), 810--833.\newline
DOI: 10.1016/j.aim.2018.03.024

\bibitem{GS04} A. Gy\'arf\'as and G. Simonyi,
\emph{Edge colorings of complete graphs without tricolored triangles},
J. Graph Theory 46 (2004), 211--216.

\bibitem{Ion09} Y. J. Ionin,
\emph{Isosceles sets},
Electron. J. Combin. 16 (2009), art. R141.\newline
DOI: 10.37236/230

\bibitem{Kid06} H. Kido,
\emph{Classification of isosceles eight-point sets in three-dimensional Euclidean space},
European J. Combin. 27 (2006), 329--341.\newline
DOI: 10.1016/j.ejc.2005.01.003

\bibitem{Kid10} H. Kido,
\emph{On isosceles sets in the $4$-dimensional Euclidean space},
Int. J. Comb. 2010 (2010), art. 803210.\newline
DOI: 10.1155/2010/803210

\bibitem{LRS77} D. G. Larman, C. A. Rogers and J. J. Seidel,
\emph{On two-distance sets in Euclidean space},
Bull. London Math. Soc. 9 (1977), 261--267.\newline
DOI: 10.1112/blms/9.3.261

\bibitem{Lis97} P. Lison\v{e}k,
\emph{New maximal two-distance sets},
J. Combin. Theory Ser. A 77 (1997), 318--338.\newline
DOI: 10.1006/jcta.1997.2749

\bibitem{Mus09} O. R. Musin,
\emph{Spherical two-distance sets},
J. Combin. Theory Ser. A 116 (2009), 988--995.\newline
DOI: 10.1016/j.jcta.2008.09.003

\bibitem{OEIS} OEIS Foundation Inc.,
\emph{The On-Line Encyclopedia of Integer Sequences},
\url{https://oeis.org}, sequence A175769.

\end{thebibliography}
\end{document}